\documentclass[11pt]{amsart}
\usepackage{amssymb}
\usepackage{amsmath}
\usepackage{color}
\newcommand{\R}{\mathbb{R}}

\newcommand{\C}{\mathbb{C}}

\font\eufm=eufm10
\def\frak#1{\hbox{\eufm#1}}
\newcommand{\bd}{\begin{document}}
\newcommand{\ed}{\end{document}}
\newcommand{\be}{\begin{enumerate}}
\newcommand{\ee}{\end{enumerate}}
\newcommand{\bi}{\begin{itemize}}
\newcommand{\ei}{\end{itemize}}
\newcommand{\ba}{\begin{array}}
\newcommand{\ea}{\end{array}}
\newcommand{\vs}{\vspace*{0.3\baselineskip}}%%%maly odstep pionowy
\newcommand{\vsm}{\vspace*{-0.3\baselineskip}}%%%maly odstep pionowy ujemny
\newtheorem{defi}{Definition}[section]
\newtheorem{tw}[defi]{Theorem}
\newtheorem{prop}[defi]{Proposition}
\newtheorem{lem}[defi]{Lemma}
\newtheorem{re}[defi]{Remark}
\newtheorem{col}[defi]{Corollary}
\newtheorem{ex}[defi]{Examples}
\newtheorem{zad}{Exercise}[section]
\newtheorem{zal}{Assumptions}[section]
\newtheorem{assumpt}[defi]{Assumptions}
\newcommand{\Om}{\Omega}
\newcommand{\om}{\omega}
\newcommand{\G}{\Gamma}
\newcommand{\D}{\Delta}
\renewcommand{\d}{\delta}
\newcommand{\ga}{\gamma}
\newcommand{\eps}{\epsilon}
\newcommand{\ove}{\overline}
\newcommand{\ms}{\oplus}
\newcommand{\mt}{\otimes}
\newcommand{\dz}{\wedge}
\newcommand{\lra}{\longrightarrow}
\newcommand{\sign}{\mbox{$ sgn $}}
\newcommand{\rel}{\mbox{$\,$\rule[0.5ex]{1.1em}{0.2pt}$\triangleright\,$}}
\newcommand{\dow}{\hspace*{\fill}\rule{1.6ex}{1.6ex}\hspace*{1em}}
\newcommand{\dowl}{\hspace*{\fill}\rule{1ex}{1ex}\hspace*{1em}}
\newcommand{\sd}{\hspace{0.3ex}\tiny{\rhd\mbox{\hspace{-2ex}}<}\hspace{0.3ex}}
\newcommand{\mmt}[2]{\mbox{$\vphantom{}_{#1}\times_{#2}$}}
\newcommand{\gotg}{\frak g}
\newcommand{\gota}{\frak a}
\newcommand{\gotb}{\frak b}
\newcommand{\gotc}{\frak c}
\newcommand{\gothh}{\frak h}
\newcommand{\gott}{\frak t}
\newcommand{\hd}{\hat{\d}}
\newcommand{\oml}{\Omega_L^{1/2}}
\newcommand{\omr}{\Omega_R^{1/2}}
\newcommand{\omh}{\Omega_c^{1/2}}
\newcommand{\lo}{\lambda_0}
\newcommand{\ro}{\rho_0}
\newcommand{\lma}{\Lambda^{max}}
\newcommand{\timh}{\times_h}
\newcommand{\Gd}{\G^{(2)}}
\newcommand{\el}{e_L}
\newcommand{\er}{e_R}
\newcommand{\GG}{\G_1\times\G_2}
\newcommand{\gdot}{\hspace{-0.1em}\cdot\hspace{-0.1em}}
\newcommand{\tran}{\frown\hspace{-2.2ex}|\hspace{1.9ex}}
\newcommand{\la}[2]{\Lambda_{#1#2}}
\newcommand{\kad}{Ad^{\#}}
\newcommand{\kkad}{\mbox{$\mathrm W$}}
\newcommand{\ad}{Ad}
\newcommand{\wl}[1]{\vphantom{X}_{#1}{\G}}
\newcommand{\te}{\tilde{e}}
\newcommand{\notka}[1]{}
\newcommand{\kom}[1]{}
\newcommand{\Mlambda}{\mbox{$\mathrm M$}}
\newcommand{\sA}{\mbox{$\mathcal A$}}
\newcommand{\sB}{\mbox{$\mathcal B$}}
\newcommand{\sC}{\mbox{$\mathcal C$}}
\newcommand{\sD}{\mbox{$\mathcal D$}}
\newcommand{\sF}{\mbox{$\mathcal F$}}
\newcommand{\sG}{\mbox{$\mathcal G$}}
\newcommand{\sL}{\mbox{$\mathcal L$}}
\newcommand{\sM}{\mbox{$\mathcal M$}}
\newcommand{\sO}{\mbox{$\mathcal O$}}
\renewcommand{\sL}{\mbox{$\mathcal L$}}
\newcommand{\sS}{\mbox{$\mathcal S$}}
\newcommand{\sT}{\mbox{$\mathcal T$}}
\newcommand{\sY}{\mbox{$\mathcal Y$}}
\newcommand{\hY}{\mbox{$\hat{Y}$}}
\newcommand{\hS}{\mbox{$\hat{S}$}}
\newcommand{\hX}{\mbox{$\hat{X}$}}
\newcommand{\dif}{differential }
\newcommand{\gru}{groupoid }
\newcommand{\grus}{groupoids }
\newcommand{\ti}{\tilde}
\newcommand{\halden}{half density }
\newcommand{\haldens}{half densities }
\renewcommand{\top}{topological }
\newcommand{\Setrel}{\mbox{\rm SetRel}}
\newcommand{\cstardwa}{\mbox{$C^*_r(\Gamma\times\Gamma)$}}
\newcommand{\cred}{\mbox{$C^*_r$}}
\newcommand{\kropka}[1]{\dot{#1}}
\newcommand{\adc}{Ad^{\gotc}}
\newcommand{\tadc}{\widetilde{Ad^{\gotc}}}
\newcommand{\adb}{Ad^{\gotb}}
\newcommand{\kot}{\mbox{ $ \Pi $}}
\newcommand{\rzutB}{\mbox{$\mathrm{P_{\gotb}}$}}
\newcommand{\rzutC}{\mbox{$\mathrm{P_{\gotc}}$}}
\newcommand{\trzutC}{\mbox{$\widetilde{\mathrm{P}}_{\gotc}$}}
\newcommand{\antyp}{\mbox{$\mathrm{S}$}}
\newcommand{\counit}{\epsilon}

\def\tgr{{\bf t}}
\def\sgr{{\bf s}}
\def\fgr{{\bf f}}
\def\rhogr{{\boldsymbol \rho}}
\begin{document}
\title{The $\kappa$-Poincar\'e Group on a  $C^*$-level.}
\author{Piotr Stachura}
\address{Faculty of Applied Informatics and Mathematics, Warsaw University of Life Sciences-SGGW,
ul Nowoursynowska 166, 02-787 Warszawa, Poland,  
e-mail: piotr\_stachura1@sggw.pl}
\date{}
\begin{abstract} The $C^*$-algebraic $\kappa$-Poincar\'{e} Group is constructed. The construction uses groupoid algebras of differential groupoids associated to Lie group decomposition.
It turns out the underlying $C^*$-algebra is the same as for  ``$\kappa$-Euclidean Group'' but a comultiplication is twisted by some unitary multiplier. Generators and commutation relations among them are presented.
\end{abstract}
\maketitle
%%%%%%%%%%%%%%%%%%%%%%%%%%%%%%%%%%%%%%%%%%%%%%%%%%%%%%%%%%%%%%%%%%%%%%%%%%%%%%%%
%%%%%%%%%%%%%%%%%%%%%%%%%%%%%%%%%%%%%%%%%%%%%%%%%%%%%%%%%%%%%%%%%%%%%%%%%%%%%%%
\section{Introduction}

The history of  $\kappa$-{\em deformation} has begun in  1992  with the work of J. Lukierski, A. Nowicki and H. Ruegg \cite{LNR}, where  this deformation 
of the enveloping algebra of Poincar\'{e} Group appeared for the first time. Next important step was the paper by S. Majid and H. Ruegg in 1994 \cite{Maj-Rueg} 
identifying bicrossedproduct structure of  $\kappa$-Poincar\'{e} algebra. 
Since then, a vast literature on the subject has been produced with many  attempts to apply this deformation to physical problems, 
but because (almost) all work was done on the level of pure algebra (or formal power series) it was hard to get more then rather formal  conclusions. This is not a review paper and we refer to 
\cite{Luk-review} (also  an  older one \cite{Bor-Pa}) for a discussion and an extensive bibliography of the subject.

This article, however,  is not about the $\kappa$-Poincar\'{e} algebra but about the $\kappa$-Poincar\'{e} Group - deformation of an algebra of functions on Poincar\'{e} Group. 
It appeared in 1994, on a  Hopf $*$-algebra level, in the work of S. Zakrzewski \cite{SZ-94}, where  it was also shown it is  a quantization of a  certain Poisson-Lie  structure. 
Soon, it became clear that this particular Poisson structure is not special for dimension 4 but has analogues in any dimension and is related to  certain 
decompositions of orthogonal Lie algebras \cite{PS-triple} and is dual to a certain Lie algebroid structure \cite{PS-poisson}. 

The main result of this work is {\em a topological version}  of $\kappa$-Poincar\'{e} Group. Since \cite{Maj-Rueg} it has been clear that had it existed it 
should have been given by some bicrossedproduct construction. The main problem is that the decomposition of a Lie algebra $\gotg=\gota\oplus\gotc$ doesn't lift to a decomposition of a Lie Group 
$G=AC$ (of course it lifts to a local decomposition, but the complement of the set of decomposable elements i.e.  $G\setminus(AC\cap CA)$ has a non empty interior) therefore the construction 
of S.~Vaes and L.~Vainerman presented in \cite{VV} can't be directly applied. 
As it has been shown already in \cite{PS-triple} there is a {\em non connected} extension of a group $A$, let's denote it by $\tilde{A}$, such that
$\tilde{A} C\cap C \tilde{A}$ is open and dense in $G$.  {\em Therefore it fits into the framework of \cite{VV} and the $\kappa$-deformation of Poincar\'{e} Group, or rather its non connected
extension, exists as a locally compact quantum group.}

Here we use approach different then used in \cite{VV}. Although less general, it has some advantages -- it is more geometric and it is easier 
to see that what we get is really a quantization of a Poisson-Lie structure. It is based on the use of groupoid algebras for differential groupoids naturally related to 
decompositions of  Lie groups.  For a global decomposition this construction was described in \cite{PS-DLG}; the result is that given a Lie group $G$ with two closed subgroups 
$B,C\subset G$ satisfying $G=BC$  one can define two differential groupoid structures on $G$ (over $B$ and $C$, this is described briefly in the second part of this introduction). 
It turns out that  $C^*$-algebras of these groupoids carry quantum group structures, in fact, sweeping under the rug some universal/reduced algebras problems, one may say that all main
ingredients of quantum group structure are just $C^*$-lifting of natural groupoid objects. As said above the situation with $\kappa$-Poincar\'{e} is not so nice, but there is  a 
global decomposition ``nearby''one can try to use; this framework was described in \cite{PSAx}.

Let us explain briefly the construction; geometric details were presented in \cite{PS-poisson}. 
By the (restricted) Poincar\'{e} Group it's  meant here the semidirect product of the (restricted)  Lorentz Group $A:=SO_0(1,n)$ and $n+1$-dimensional vector Minkowski space. It turns out, 
that it can be realized as a subgroup $(TA)^0$ of $T^*G$ for  $G:=SO_0(1,n+1)$, 
where $A$ is embedded naturally into $G$ (as the stabilizer of a spacelike vector). On $(TA)^0$ there is a Poisson structure
dual to a Lie algebroid structure related to the decomposition of $\gotg$ -- the Lie algebra of $G$ into two subalgebras $\gotg=\gota\oplus \gotc$, 
where $\gota$ is the Lie algebra of $A$.  This Lie algebroid is the algebroid of the Lie groupoid $AC\cap CA\subset G$, here $C$ is the Lie subgroup with algebra $\gotc$. As said above, 
the set $AC\cap CA$ is too small, but one can find $\tilde{A}\subset G$ -- a non-connected extension of $A$, in fact this is the normalizer of $A$ in $G$, 
such that  $\Gamma:=\tilde{A} C\cap C \tilde{A}$ is open and dense in $G$. The set $\Gamma$ is a differential groupoid (over $\tilde{A}$) and $C^*$-algebra of this groupoid is 
the $C^*$-algebra of $\kappa$-Poincar\'{e} Group. The method used in \cite{PSAx} relies essentially on the fact that the algebra $\gotc$ has a second complementary algebra $\gotb$ and
the decomposition $\gotg=\gotc\oplus\gotb$ lifts to {\em a global decomposition} $G=BC$ and this is just the Iwasawa decomposition (i.e. $B=SO(n+1)$). 
This global decomposition defines a groupoid  $G_B:G\rightrightarrows B$ and its $C^*$-algebra is the underlying algebra of the quantum group which may be 
called ``Quantum $\kappa$-Euclidean Group''.  Our groupoid $\Gamma$ embeds into $G_B$ in such a way that it is possible to prove that their $C^*$-algebras are 
the same but the comultiplication of $\kappa$-Poincar\'{e} is 
comultiplication of $\kappa$-Euclidean twisted by a  unitary multiplier. So one may say that  {\em quantum spaces underlying $\kappa$-Poincar\'{e} and $\kappa$-Euclidean 
groups are the same and only group structures are different.}

The embedding $\Gamma\hookrightarrow G_B$ essentially is given by embedding of $\tilde{A} \hookrightarrow SO(n+1)$ as a dense open subset, so it is a kind of compactification of 
$\tilde{A}$ (which consists of two copies of the (restricted) Lorentz Group $SO_0(1,n)$). This compactification solves the problem, that some natural operators, that 
``should be'' self-adjoint elements affiliated with $C^*(\Gamma)$ are defined by non complete vector fields, so they are not essentially self-adjoint on their ``natural'' domains and 
this embedding just defines ``correct'' domains.

In this work we consider only $C^*$-algebra with comultiplication, and not discuss other ingredients like antipode, Haar weight, etc.; they can be constructed using methods presented in 
\cite{PS-DLG}.

In the remaining part of the Introduction we recall basics of  groupoid algebras, groupoids related to decomposition of groups and results of \cite{PSAx} essential in the following.
The short description of the content of each section  is given at the end of the Introduction.
%%%%%%%%%%%%%%%%%%%%%%%%%%%%%%%%%%%%%%%%%%%%%%%%%%%%%%%%%%%%%%%%%%%%%%%%%%%%%%%%%%%%%%%%%%%%%%%%%%%%%%%%%%%%%%%
\subsection{Groupoid algebras} We will use groupoid algebras, so now we recall basic facts and establish the 
relevant notation. We refer to \cite{PS-DLG,DG} for a detailed exposition and to \cite{PSAx} for basics of formalism. 

All {\em manifolds} are smooth, Hausdorff, second countable and {\em submanifolds} are embedded.
For a manifold $M$  by $\omh(M)$ we denote the vector space of smooth, compactly supported, complex \haldens  on $M$; it is equipped with the  scalar product
$\displaystyle (\psi_1\,|\,\psi_2):=\int_M\overline{\psi_1}\psi_2$ and $L^2(M)$ is the completion of $\omh(M)$ in the associated norm. Clearly, if we choose
some $\psi_0$ -- non vanishing, real \halden on $M$,  there is the equality $\omh(M)=\{f\,\psi_0\,,\,f\in \sD(M)\}$, 
where $\sD(M)$ stands for smooth, complex and compactly supported functions on $M$.

Let $\Gamma\rightrightarrows E$ be a differential groupoid. By $e_R$ ($e_L$) we denote  the {\em source (target) projection} 
 and called it {\em right (left)  projection}; a groupoid inverse is denoted by $s$.  Let $\Om_L^{1/2}, \Om_R^{1/2}$ denote bundles of complex  \haldens 
 along left and right fibers. {\em A groupoid *-algebra} $\sA(\Gamma)$ is a vector space of compactly 
supported, smooth sections of $\Om_L^{1/2}\mt\Om_R^{1/2}$  together with a convolution and $*$-operation.
To write explicit  formulae let us choose $\lambda_0$ - a real, non vanishing, left invariant \halden along left fibers 
(in fact that means we choose a Haar system on $\Gamma$, however nothing depends on this choice, details are  given in \cite{DG}.) 
Let $\rho_0:=s(\lambda_0)$ be the corresponding right invariant \halden
and $\om_0:=\lambda_0\mt \rho_0$. Then any $\om\in\sA(\Gamma)$ can be written as $\om=f\om_0$ for a  unique function 
$f\in\sD(\Gamma)$. With such a choice  we write $(f_1\om_0)\,(f_2\om_0)=:(f_1*f_2)\om_0$,
$(f\om_0)^*=:(f^*)\om_0$ and: \notka{group-alg-mult}
%%%
\begin{equation}\label{group-alg-mult}
(f_1*f_2)(\gamma):=\int_{F_l(\gamma)} \lambda_0^2(\gamma') f_1(\gamma')f_2(s(\gamma') \gamma)=
\int_{F_r(\gamma)} \rho_0^2(\gamma')f_1(\gamma s(\gamma')f_2(\gamma')\,\,,\,\,f^*(\gamma):=\overline{f(s(\gamma))}
\end{equation}
$F_l(\gamma)$ and $F_r(\gamma)$ are left and right fibers passing through $\gamma$ e.g.  $F_l(\gamma):=e_L^{-1}(e_L(\gamma))$. \\
The choice of $\om_0$ defines a norm that makes  $\sA(\Gamma)$ a normed $*$-algebra:
$$||f\om_0||_0=:||f||_0=max\left\{\sup_{e\in E}\int_{\el^{-1}(e)}\lambda_0^2(\gamma)|f(\gamma)|,\,
\sup_{e\in E}\int_{\er^{-1}(e)}\rho_0^2(\gamma)|f(\gamma)|\right\}$$

There is a faithful representation $\pi_{id}$ of $\sA(\Gamma)$ on $L^2(\Gamma)$ described as follows: choose $\nu_0$ - a real, 
non vanishing \halden  on $E$; since $\er$ is a surjective submersion one can define $\psi_0:=\rho_0\mt\nu_0$ - this is 
a real, non vanishing, \halden on $\Gamma$. For $\psi=f_2\psi_0,\, f_2\in\sD(\Gamma)$  the representation is given by 
$\pi_{id}(f_1\om_0) (f_2\psi_0)=:(\pi_{id}(f_1)f_2)\psi_0$ and $\pi_{id}(f_1)f_2=f_1*f_2$  is as in (\ref{group-alg-mult}). 
The estimate $||\pi_{id}(\om)||\leq||\om||_0$ makes  possible the definition:
{\em  The reduced $C^*$-algebra of a groupoid} -- $C^*_r(\Gamma)$ is the completion of $\sA(\Gamma)$ in the  norm $||\om||:=||\pi_{id}(\om)||$.
We will also use the following fact which is a direct consequence of the definition of the norm $||f||_0$:
\begin{lem}\label{lemma-ind-lim}\notka{lemma-ind-lim} Let $U\subset \Gamma$ be an open set with compact closure. There exists $M$ such that
$||f||_0\leq M \sup |f(\gamma)|$ for any $f\in\sD(\Gamma)$ with support in $U$. If $f_n\in\sD(\Gamma) $ 
have supports in a fixed compact set and $f_n$ converges to $f\in\sD(\Gamma)$ uniformly then $f_n\om_0$ converges 
to $f\om_0$ in  $C^*_r(\Gamma)$.
\dowl
\end{lem}
%%%%%%%%%%%%%%%%
We will use {\em Zakrzewski category} of groupoids \cite{SZ1,SZ2}. {\em Morphisms} are not mappings  (functors)  but relations satisfying certain natural properties. 
A morphism $h:\Gamma\rel\Gamma'$ of differential groupoids defines a mapping 
$\hat{h}:\sA(\Gamma)\rightarrow L(\sA(\Gamma'))$ (linear mappings of $\sA(\Gamma')$), 
this mapping commutes with (right) multiplication in $\sA(\Gamma')$ i.e. 
$$\hat{h}(\om)(\om')\om''=\hat{h}(\om)(\om' \om'')\,\,,\,\om\in\sA(\Gamma)\,,\,\om',\om''\in\sA(\Gamma')$$
and we use  notation $\hat{h}(\om)\om'$ (see formula (\ref{app-hat-delta0-form}) in the Appendix for an example of such a mapping);
 there is also a representation $\pi_h$ of $\sA(\Gamma)$ on 
$L^2(\Gamma')$; these objects satisfy some obvious compatibility conditions with respect to multiplication 
and $*$-operation (see \cite{DG} for details).

{\em A bisection}  $B$ of $\Gamma\rightrightarrows E$ is a submanifold such that $e_L|_B, e_R|_B :B\rightarrow E$ are diffeomorphisms. Bisections act on $\cred(\Gamma)$  as unitary 
multipliers and are transported by morphisms:  if $B\subset \Gamma$ is a bisection and $h:\Gamma\rel\Gamma'$ is a morphism then the set $h(B)$ is a bisection of $\Gamma'$.

\subsection{Group decompositions, related groupoids and quantum groups} 
Now we briefly recall some facts about {\em double groups}. Let $G$ be a group and $A,B\subset G$ subgroups 
such that $A\cap B=\{e\}$.  Every element $g$ in the set   $\Gamma:=AB\cap BA$ can be written uniquely as 
$$g=a_L(g) b_R(g)=b_L(g) a_R(g)\,,\,a_L(g),  a_R(g)\in A\,,\,b_L(g),b_R(g)\in B.$$
These decompostions define surjections: $a_L,a_R: \Gamma\rightarrow A$ and 
$b_L,b_R: \Gamma\rightarrow B$  (in fact $a_L, b_R$ are defined on $AB$ and $b_L, a_R$ on $BA$, 
we will denote these extensions by the same symbols). The formulae:
\begin{align*}E &:=A \,, \quad s(g):= b_L(g)^{-1}a_L(g)=a_R(g) b_R(g)^{-1}\,,\,\\
Gr(m) &:= \{(b_1 a b_2; b_1 a, a b_2) :b_1a, ab_2\in\Gamma\}
\end{align*}
define the structure of the  groupoid $\Gamma_A:\Gamma\rightrightarrows A$; 
the analogous formulae define the groupoid  $\Gamma_B: \Gamma\rightrightarrows B$. 
On the other hand for a subgroup $B\subset G$ there is a (right) transformation groupoid 
$(B\setminus G)\rtimes B$. The following lemma  \cite{PSAx} explains relation between these groupoids.
\begin{lem}\label{lemma-embed}\notka{lemma-embed} The map:
$$\Gamma_A\ni g\mapsto ([a_L(g)], b_R(g))\in (B\backslash G)\rtimes B$$
is an isomorphism of the groupoid $\Gamma_A$ with the restriction of a 
(right) transformation groupoid $ (B\backslash G)\rtimes B$ to the set $\{[a] : a\in A\}\subset B\backslash G$.
\dowl
\end{lem}

If $AB=G$ (i.e. $\Gamma=G$)  the  triple $(G;A,B)$ is called {\em a double group} 
and in this situation we will denote groupoids 
$\Gamma_A, \Gamma_B$ by $G_A, G_B$. It turns out that {\em the transposition of multiplication relation $m_B$} i.e. $\delta_0:=m_B^T: G_A\rel G_A\times G_A$ is a 
coassociative morphism of groupoids. 
Applying the lemma \ref{lemma-embed} to the groupoid $G_A$ we can identify it with the transformation groupoid 
$(B\backslash G)\rtimes B$. So $G_A=A\rtimes B$ is a right transformation groupoid for the action
$(a,b)\mapsto a_R(ab)$ i.e the structure is given by:
$$E:=\{(a,e): a\in A\}\,,\,s(a,b):=(a_R(ab),b^{-1}),$$
$$m:=\{(a_1,b_1 b_2;a_1, b_1,a_R(a_1 b_1),b_2): a_1\in A\,,\,b_1,b_2\in B\}$$
{\em In the formula above, we identified a relation $m:\Gamma\times \Gamma\rel\Gamma$ with its graph, i.e. subset of 
 $\Gamma\times \Gamma\times \Gamma$. We will use such notation throughout the paper.}
If $G$ is a Lie group, $A,B$ are closed subgroups, $A\cap B=\{e\}\,,\,AB=G$ 
then $(G;A,B)$ is called {\em a double Lie group}, abbreviated in the following as DLG.
It turns out that the mapping $\widehat{\delta_0}$, defined by the morphism $\delta_0$ (compare (\ref{app-hat-delta0-form}) in Appendix),
extends to the coassociative morphism $\Delta$ of $C^*_r(G_A)$ and $C^*_r(G_A\times G_A)=C^*_r(G_A)\mt C^*_r(G_A)$ 
which satisfies {\em density conditions}:
\begin{equation*}
cls\{\Delta(a)(I\mt b) :a,b \in C^*_r(G_A)\}=cls\{\Delta(a)(b\mt I) :a,b \in C^*_r(G_A)\}=C^*_r(G_A)\mt C^*_r(G_A),
\end{equation*}
where {\em cls} denotes the closed linear span. 
There are other objects that make the pair $(C^*_r(G_A), \Delta)$ a locally compact quantum group; we refer to  \cite{PS-DLG} for details.

\subsection{Framework for $\kappa$-Poincar\'{e}}
The framework for $\kappa$-Poincare group we are going to use was presented  in \cite{PSAx}. Let us now recall basic facts established there.
Let $G$ be a group and $A,B,C\subset G$ subgroups satisfying conditions: 
\notka{rozklad-warunki}\begin{equation*}\label{rozklad-warunki} 
B\cap C=\{e\}=A\cap C\,,\,B C=G.
\end{equation*}
i.e. $(G;B,C)$ is a double group. As described above, in this situation, there is the groupoid $G_B$, and the (coassociative) morphism
$\delta_0:G_B\rel G_B\times G_B$; explicitly the graph of $\delta_0$ is equal to: 
\notka{def-delta0}\begin{equation}\label{def-delta0} \delta_0=\{(b_1 c, c b_2;b_1 c b_2)
  \,:\, b_1,b_2\in B\,,\,c\in C\}
\end{equation}
%%%%
Using the lemma \ref{lemma-embed} we see that this is a transformation groupoid $(C\backslash G) \rtimes C$ and the isomorphism is
$$(C\backslash G) \times C\ni ([g],c)\mapsto b_R(g)c\in G$$
%%%%%%%%
Let $\Gamma:=A C\cap C A$ and  consider on $\Gamma$ the groupoid structure $\Gamma_A: \Gamma\rightrightarrows A$ described above, together with a relation
(the transposition of the multiplication in $\Gamma_C:\Gamma\rightrightarrows C$):
\begin{equation*}\label{def-tmct} \notka{def-tmct} \tilde{m}_C^T:=\{(a_1 c_1, c_1 a_2 ; a_1 c_1 a_2) : a_1 c_1, c_1 a_2\in\Gamma\}\subset \Gamma\times\Gamma\times\Gamma.
\end{equation*}
The corresponding projections will be denoted by $\tilde{c}_L, \tilde{c}_R $ and 
$a_R, a_L$. Again, by the lemma \ref{lemma-embed} we identify the groupoid $\Gamma_A$ with the restriction 
of $(C\backslash G) \rtimes C$ and then with the restriction
of $G_B$  to the set $B':=B\cap CA$, i.e. with $b_L^{-1}(B')\cap b_R^{-1}(B')$. 
This restriction will be denoted by $\Gamma_{B'}$ (instead of more adequate but rather inconvenient $G_B|_{B'}$). 
This isomorphism and its inverse are given by:\notka{wzor-embed}
\begin{align}\label{wzor-embed}
\Gamma_A\ni & a c\mapsto b_R(a) c\in \Gamma_{B'}\quad\quad,\quad \Gamma_{B'}\ni b c\mapsto  a_R(b) c\in \G_A
\end{align}
The image of $\tilde{m}_C^T$ inside $\Gamma_{B'}\times \Gamma_{B'}\times \Gamma_{B'}$ is equal to:
\begin{equation*}\label{mct}
\{(b_R(a_1) c_1, b_R(a_2) c_2; b_R(a_1 a_2) c_2): a_1 c_1=\tilde{c}_1 \tilde{a}_1, c_1\tilde{a}_2=a_2 c_2\}
\end{equation*}
%%%%%%%%%%%%%%%%%%%%%%%%%%%%%%%%%%
The following object plays the major role in what follows:
\notka{def-twist}\begin{equation}\label{def-twist}
T:=\{(g,b):c_R(g)b\in A\}=\{(b_1 \tilde{c}_L(b_2)^{-1}, b_2) : b_1\in B,b_2\in B'\}\subset G_B\times G_B.
\end{equation}
Using the definition (\ref{def-delta0}) of $\delta_0$   one easily computes images of $T$ by 
relations $id\times\delta_0$ and $\delta_0\times id$:
%%%%%%%%%%%%%%%%%%%%%%
\begin{align*}
(id\times \delta_0)& T=\{(g_1,b_2,b_3): b_2 b_3\in B'\,,\,c_R(g_1)=\tilde{c}_L(b_2 b_3)^{-1}\}\\
(\delta_0\times id)T& =\{(g_1,g_2,b_3): c_R(g_1)=c_L(g_2)\,,\,b_3\in B'\,,\,c_R(g_2)=\tilde{c}_L(b_3)^{-1}\}
\end{align*}
%%%%%%%%%%%%%%%%%%%%%%%%%%%%%%
Let us also denote $T_{12}:=T\times B\subset G_B\times G_B\times G_B$ and 
$T_{23}:=B\times T\subset G_B\times G_B\times G_B$. \\
%%%%%%%%%%%%%%
Main properties of $T$ are listed in the following lemma (proven in \cite{PSAx}):
\begin{prop} \label{prop-twist}\notka{prop-twist}\begin{enumerate}
\item $T$ is a section of left and right projections (in $G_B\times G_B$) over the set $B\times B'$ and a 
bisection of $G_B\times \Gamma_{B'}$;
\item $(id\times \delta_0) T$ is a section of left and right projections 
(in $G_B\times G_B\times G_B$) over the set 
$B\times \delta_0(B')=\{(b_1,b_2,b_3) : b_2 b_3\in B'\}$;
\item $(\delta_0\times id) T$ is a section of left and right projections  over the set 
$B\times B \times B'$;
\item   $T_{23}(id\times \delta_0)T= T_{12}(\delta_0\times id) T$ (equality of sets in 
$G_B\times G_B\times G_B$), 
moreover this set is a section of the right projection over $B\times (\delta_0(B')\cap(B\times B'))$ and 
the left projection over
$B\times B'\times B'$.
\end{enumerate}
\end{prop}
\noindent
Due to this proposition the left multiplication by $T$, which we denote by the same symbol, is a 
bijection of $G_B \times b_L^{-1}(B')$.
Let $Ad_{T}:G_B\times G_B\rel G_B\times G_B$ be a relation defined by:\notka{def-AdT}
\begin{equation}\label{def-AdT}
(g_1,g_2;g_3,g_4)\in Ad_{T}\iff \exists t_1, t_2\in T : (g_1,g_2)=t_1(g_3,g_4)(s_B\times s_B)(t_2).
\end{equation}
and let us define  the relation $\delta:=Ad_T\cdot \delta_0: G_B\rel G_B\times G_B$
\notka{def-delta}\begin{equation}\label{def-delta} 
\begin{split}
\delta= \{ & (b_R(b_3 b_2^{-1} \tilde{c}_L(b_2))\tilde{c}_L(b_2)^{-1} c_L(b_2 c_2) \tilde{c}_L(b_R(b_2 c_2)), 
b_2 c_2; b_3 c_2) :\\
 & b_3\in B, c_2\in C,  b_2, b_R(b_2 c_2)\in B' \}.
\end{split}
\end{equation}
The relation between $\delta$ and $\tilde{m}_C^T$ is explained in the lemma: \cite{PSAx}
\begin{lem} $\delta$ is an extension of $\tilde{m}_C^T$ i.e. $\tilde{m}_C^T\subset \delta$
\end{lem}
%%%%%%

\noindent
Addition of some differential conditions to this situation makes possible to use $T$ to twist the comultiplication on $C^*_r(G_B)$: 

\begin{assumpt}\label{basic-assumpt}\notka{basic-assumpt}
\begin{enumerate}
\item $G$ is a Lie group and $A,B,C$ are closed Lie subgroups such that \\
$B\cap C=\{e\}=A\cap C\,,\,B C=G$ (i.e. $(G;B,C)$ is a DLG).
\item The set $\Gamma:=C A\cap AC$ is open and dense in $G$.
\item Let $U:=b_L^{-1}(B')$ and $\sA(U)$ be the linear space of elements from $\sA(G_B)$ supported in $U$.
We assume that $\sA(U)$ is dense in $C^*_r(G_B)$.
\item For a compact set $K_C\subset C$, open $V\subset B$ and $(b_1, b_2)\in B\times B'$ let us define a set
$Z(b_1,b_2,K_C;V):=K_C\cap\{c\in C: b_R(b_1 c) b_2\in V\}$ and a function:
$$B\times B' \ni(b_1,b_2)\mapsto \mu(b_1,b_2,K_C;V):=\int_{Z(b_1,b_2,K_C;V)} d_l c.$$
For compact sets $K_1\subset B$ and $K_2\subset B'$ let $\mu(K_1,K_2,K_C;V):=\sup\{\mu(b_1,b_2,K_C;V): b_1\in K_1\,,\,b_2\in K_2\}$
We assume that \begin{equation*}
\forall\,\epsilon>0\,\exists\, V-{\rm a\, neighborhood\,of\,}B\setminus B'{\rm\, in\, B}\,:\, \mu(K_1,K_2,K_C;V)\leq\epsilon
\end{equation*}
\end{enumerate}
\end{assumpt}
\begin{re}\label{remark-B}\notka{remark-B}
\begin{itemize}
\item It follows from the first and the second assumptions that $B'$ is open and dense in $B$. 
\item The second assumption can be replaced by the following  two conditions:
\begin{itemize}
\item[a)] $\gotg=\gota\oplus\gotc$, where $\gotg, \gota,\gotc$ are lie algebras of $G,A,C$, respectively. (then $AC$ and $CA$ are open);
\item[b)] $AC\cap CA$ is dense in $G$.
\end{itemize}
\item These assuumptions are not very pleasent and probably, at least some of them, redundant; but they are sufficient to get results in case of quantum ``ax+b''  
group and $\kappa$-Poincar\'{e}. This framework, however, doesn't work for dual to $\kappa$-Poncar\'{e} and at the moment it is unclear if 
and how that example may be handled (with this approach).
\end{itemize}
\end{re}

The following proposition was proven in \cite{PSAx}:
%%%%%%%%%%%%%
\begin{prop}\label{prop-Delta} \notka{prop-Delta} Assume the conditions listed in (\ref{basic-assumpt}) are satisfied. Then
\begin{itemize}
\item[a)] $C^*_r(\Gamma_{B'})=C^*_r(G_B)$ -- for this equality,  it is sufficient to satisfy (1),(2),(3) from \ref{basic-assumpt};
\item[b)] The mapping
  $T:\sA(G_B\times U)\rightarrow\sA(G_B\times U)$
extends to the unitary $\widehat{\sT}\in M(C^*_r(G_B)\mt C^*_r(G_B))$ which satisfies:
$$(\widehat{\sT}\mt I)(\Delta_0\mt id)\widehat{\sT}=(I \mt \widehat{\sT})(id\mt \Delta_0)\widehat{\sT}$$
\item[c)] Because of b), the formula $\Delta(a):=\widehat{\sT}\Delta_0(a)\widehat{\sT}^{-1}$ defines a coassociative morphism. 
For this morphism ("cls" stands for "closed linear span"):\notka{density-c}
\begin{equation} \label{density-c}
cls\{\Delta(a)(I\mt c)\,,\,a,c\in C^*_r(G_B)\}=cls\{\Delta(a)(c\mt I)\,,\,a,c\in C^*_r(G_B)\}=
C^*_r(G_B)\mt C^*_r(G_B).
\end{equation}
\end{itemize}
\dow
\end{prop}

The second  section describes the situation for $\kappa$-Poincar\'{e}, i.e.~we define  groups $G,A,B,C$,  compute explicit formulae for decompositions and describe 
structure of the  groupoid $G_B$. In the third one,  we verify Assumptions \ref{basic-assumpt} and, by the Prop.~\ref{prop-Delta} get the $C^*$-algebra and 
comultiplication for $\kappa$-Poincar\'{e} Group. The fourth section describes generators of this $C^*$-algebra, 
computes commutation relations and comultiplication on generators; also, the twist is described in more details. 
In the last but one section  we compare our formulae to  the ones in  \cite{SZ-94} and in the last one we discuss ``quantum $\kappa$-Minkowski Space''.
Finally there is an  Appendix with some formulae needed here and  proven elsewhere.
%%%%%%%%%%%%%%%%%%%%%%%%%%%%%%%%%%%%%%%%%%%%%%%%%%%%%%%%%%%%%%%%%%%%%%%%%%%%%%%%%%%%%%%%%%%%%%%%%%%%%%%%%%%%%%%%%%%%%%%%%%%%%%%%%%%%%%%%%%%%%5
%%%%%%%%%%%%%%%%%%%%%%%%%%%%%%%%%%%%%%%%%%%%%%%%%%%%%%%%%%%%%%%%%%%%%%%%%%%%%%%%%%%%%%%%%%%%%%%%%%%%%%%%%%%%%%%%%%%%%%%%%%%%%%%%%%%%%%%%%%%%%
%%%%%%%%%%%%%%%%%%%%%%%%%%%%%%%%%%%%%%%%%%%%%%%%%%%%%%%%%%%%%%%%%%%%%%%%%%%%%%%%%%%%%%%%%%%%%%%%%%%%%%%%%%%%%%%%%%%%%%%%%%%%%%%%%%%%%%%%%%%%
\section{Decompositions defining $\kappa$-Poincar\'{e} Group }\label{sec:decomp}

Let $(V,\eta)$ be $n+2, n\geq1$  dimensional vector Minkowski space (signature is $(+,-,\dots,-)$).
Let us choose an orthonormal basis $(e_0, e_1,\dots,e_{n+1})$ 
and identify the (special) orthogonal group $SO(\eta)$ with the corresponding group of matricies $SO(1,n+1)$.

Let $G:=SO_0(1,n+1)\subset SO(1,n+1)$ be the connected component of identity and  $\gotg:=span\{\Mlambda_{\alpha \beta}\,,\,\alpha,\beta=0,\dots, n+1\}$ 
be its Lie algebra (see Appendix for notation). Consider three closed subgroups $A,B,C\subset G$:

{\em The group $A$} is a non-connected extension of $SO_0(1,n)$ inside $SO_0(1,n+1)$ defined by: \notka{def-A}
\begin{equation}\label{def-A}
A:=\left\{\left(\begin{array}{cc} 
u & 0\\
0 & 1
\end{array}\right)
\left(\begin{array}{cc} 
I_n &0\\
0 & h
\end{array}\right)\,:\,\,u\in SO_0(1,n),\,h:=\left(\begin{array}{cc} 
d &0\\
0 &  d
\end{array}\right)\,,\,d=\pm 1 \right\}.
\end{equation}
In fact, it is not hard to see that $A$ is the normalizer of $SO_0(1,n)$ (embedded into upper left corner) inside $SO_0(1,n+1)$. 
The Lie algebra of $A$ is $\gota:=span\{\Mlambda_{0 m}\,,\,m=1,\dots, n\}$.

\noindent
We parameterize $SO_0(1,n)$ by: 
$$\{z\in \R^n: |z|<1\}\times SO(n)\ni (z,U)\mapsto \left(\begin{array}{cc}
\frac{1+|z|^2}{1-|z|^2} & \frac{2}{1-|z|^2}z^tU\\ 
\frac{2 }{1-|z|^2} z & (I+\frac{2}{1-|z|^2}z z^t)U\end{array}\right)\in
SO_0(1,n)$$
\begin{re}
This is the standard  $boost\times rotation$ parametrization of the Lorentz Group; the parameter $z$ is related to 
a velocity $v$ by:
$$z=\frac{v}{1+\sqrt{1-|v|^2}}\,\,,\,\,\,v=\frac{2 z}{1+|z|^2}$$
\notka{jezeli $|v|=\tanh \chi$, to $|z|=\tanh(\frac\chi2)$}
\end{re}
\noindent
In this way we obtain the parametrization of $A$: \notka{def-zUd}
\begin{equation}\label{def-zUd}
\{z\in \R^n: |z|<1\}\times SO(n)\times\{-1,1\}\ni(z,U,d)\mapsto  \left(\begin{array}{ccc}
\frac{1+|z|^2}{1-|z|^2} & \frac{2}{1-|z|^2}z^tU D & 0\\ 
\frac{2 }{1-|z|^2} z & (I+\frac{2}{1-|z|^2}z z^t)U D & 0\\
0 & 0 & d\end{array}\right), \end{equation}
where $D=\left(\begin{array}{cc}I_{n-1} & 0\\0 & d\end{array}\right)$; we will also denote by $(z,U,d)$ the corresponding element of $G$.

{\em The group $B$} is $SO(n+1)$ embedded into $G$ by: $SO(n+1)\ni g\mapsto \left(\begin{array}{cc} 1 & 0 \\0 & g\end{array}\right)\in G$; its Lie algebra is
$\gotb:=\{\Mlambda_{kl}\,,\,k,l=1,\dots, n+1\}$.
Elements of $SO(n+1)$ will be written as: \notka{def-B}
\begin{equation}\label{def-B}
(\Lambda,u,w,\alpha):=\left(\begin{array}{cc} 
\Lambda  & u\\
 w^t & \alpha
\end{array}\right),\quad\Lambda\in M_n(\R) , u,w\in\R^n\,,\,\alpha\in[-1,1],
\end{equation}
and $\Lambda, u, w, \alpha$ satisfy:  \notka{param-SOn}
\begin{equation}\label{param-SOn}\Lambda \Lambda^t+u u^t=I\,,\,\Lambda w+\alpha u=0\,,\,\Lambda^t \Lambda+ w w^t=I\,,\,
\Lambda^t u+\alpha w=0\,,\,|u|^2+\alpha^2=|w|^2+\alpha^2=1;
\end{equation}
these  equations imply that $\alpha=\det(\Lambda)$. Again we will denote by $(\Lambda,u,w,\alpha)$ the corresponding element of $G$.

{\em The group $C$}  is: \notka{def-C}
\begin{equation}\label{def-C}
C:=\left\{\left(\begin{array}{ccc} 
\frac{s^2+1+|y|^2}{2 s} & -\frac1s y^t & \frac{s^2-1+|y|^2}{2 s}\\
-y & I & -y\\
\frac{s^2-1-|y|^2}{2 s} &  \frac1s y^t & \frac{s^2+1-|y|^2}{2 s}
\end{array}\right)\,s\in\R_+, y\in\R^n\right\}\subset G\,;
\end{equation}
it is isomorphic to the semidirect product of $\R_+$ and $\R^n$  $\,\{(s,y)\in \R_+\times \R^n\}$  with  multiplication
$\displaystyle (s_1,y_1) (s_2,y_2):=(s_1 s_2, s_2 y_1+y_2)$.
As before  we will use $(s,y)$ to denote the corresponding element of $G$. The Lie algebra of $C$ is 
$\gotc:=span\{ \Mlambda_{\beta 0}-\Mlambda_{\beta (n+1)}\,,\,\beta=0,\dots,  n\}=span\{ \Mlambda_{k 0}-\Mlambda_{k(n+1)}\,,\,k=1,\dots, n+1\}$.
The coordinates $(s,y)$ are related to basis in $\gotc$ as: \notka{def-wsp-sy}
\begin{align}\label{def-wsp-sy}
(s,y) &=\exp(-(\log s)  \Mlambda_{0(n+1)}) \exp(\Mlambda(y))\,,& \Mlambda(y)&:=\sum_{k=1}^n y_k(\Mlambda_{k(n+1)} -\Mlambda_{k 0})
\end{align}
\begin{re}
More geometric description of data defining groups $A,B,C$ was given in \cite{PS-poisson}; 
essentially we have to choose two orthogonal  vectors in $n+2$ dimensional Minkowski(vector) space, one spacelike and one timelike.
\end{re}
%%%%%%%%%%%%%%%%%%%%%%%%%%%%%%%%%%%%%%%%%%%%%%%%%%%%%
The Iwasawa decomposition for $G$ is $G=B C=C B$. In the following we will need explicit relation between two forms of this decomposition
 i.e. solutions of the equation \notka{eq-iwasawa}
\begin{align}\label{eq-iwasawa} 
(\Lambda,u,w,\alpha)(s,y)&=(\tilde{s},\tilde{y})(\tilde{\Lambda}, \tilde{u}, \tilde{w}, \tilde{\alpha})
\end{align}
\begin{lem}\label{lemat-iwasawa}\notka{lemat-iwasawa}
(a) Let  $(s,y)\in C$ and $(\Lambda,u,w,\alpha)\in B$.  The equation (\ref{eq-iwasawa}) has the (unique) solution 
$(\tilde{s},\tilde{y})\in C$ and $(\tilde{\Lambda}, \tilde{u}, \tilde{w}, \tilde{\alpha})\in B$ given by formulae:
\notka{solution-iwasawa1}
%%%%%%%%%%%%%%%%%%%%%%%%%%%%%%%%%%%%%%%%%%%
\begin{align}
\tilde{s} &=M s =-w^ty+\frac{s^2+1+|y|^2}{2 s}+\alpha\frac{s^2-1-|y|^2}{2 s}  &  \tilde{y} &=\Lambda y-\frac{s^2-1-|y|^2}{2 s} u\nonumber\\
\tilde{u} &=\frac{1}{M s}\left((1-\frac{w^t y}{s}) u -\frac{1-\alpha}{s} \Lambda y\right)  &  \tilde{w} &=\frac{1}{M s}\left(w-\frac{1-\alpha}{s}  y\right)\nonumber
\end{align}
\begin{align}\label{solution-iwasawa1}
\tilde{\Lambda} &=\left[\Lambda-\frac1{\alpha-1}u w^t\right]\left[I-\frac1{M (1-\alpha)}(w-\frac{1-\alpha}{s}  y)(w^t-\frac{1-\alpha}{s}  y^t)\right] \\
\tilde{\alpha} &=1-\frac{1-\alpha}{M s^2}\,,\,\, {\rm where}\nonumber
\end{align}\vspace{-4ex}
\begin{align*}
\begin{split}
M & :=\frac1{2(1-\alpha)} \left(\left(\frac{1-\alpha}{s}\right)^2+ |w-\frac{1-\alpha}{s} y|^2\right)=\\
& = \frac{1}{2}\left(\frac{1}{s^2}+\frac{|y|^2}{s^2}+1\right)-\frac{\alpha}{2}\left(\frac{1}{s^2}+\frac{|y|^2}{s^2}-1\right) -\frac{w^t y}{s}
\end{split}
\end{align*}
%%%
(b) Let $(\tilde{s},\tilde{y})\in C$ and  $(\tilde{\Lambda}, \tilde{u}, \tilde{w}, \tilde{\alpha})\in B$. The equation (\ref{eq-iwasawa}) has  the (unique) solution 
$(s,y)\in C$ and  $(\Lambda,u,w,\alpha)\in B$ given by formulae:
\notka{solution-iwasawa2} 
\begin{align}\label{solution-iwasawa2}
s&=\frac{\tilde{s}}{\tilde{M}}=
\frac{2 \tilde{s} (1-\tilde{\alpha})}{\tilde{s}^2 (1-\tilde{\alpha})^2+|\tilde{u}+(1-\tilde{\alpha})\tilde{y}|^2} &
y&=\frac{1}{\tilde{M}}\left(\tilde{\Lambda}^t \tilde{y}-\frac{\tilde{s}^2+|\tilde{y}|^2-1}{2} \tilde{w}\right)\nonumber\\
u&=\frac{\tilde{s}}{\tilde{M}}\left(\tilde{u}+(1-\tilde{\alpha})\tilde{y}\right) &
w&=\frac{\tilde{s}}{\tilde{M}}\left((1-\tilde{\alpha})\tilde{\Lambda}^t\tilde{y}+(1+\tilde{u}^t\tilde{y}) \tilde{w}\right)
\end{align}
\begin{align}
\Lambda&=\tilde{\Lambda}-\frac{1}{\tilde{M}}\left[(1+\tilde{u}^t\tilde{y})\tilde{y}\tilde{w}^t-
\frac{\tilde{s}^2+|\tilde{y}|^2-1}{2} \tilde{u}\tilde{w}^t+(\tilde{u}+(1-\tilde{\alpha})\tilde{y})\tilde{y}^t\tilde{\Lambda}\right]=  \nonumber\\
&=\left[I-\frac{1}{\tilde{M}(1-\tilde{\alpha})}(\tilde{u}+(1-\tilde{\alpha})\tilde{y})(\tilde{u}+(1-\tilde{\alpha})\tilde{y})^t\right]
\left[\tilde{\Lambda}+\frac{1}{1-\tilde{\alpha}}\tilde{u}\tilde{w}^t\right]\nonumber\\
\alpha& =1-\frac{\tilde{s}^2(1-\tilde{\alpha})}{\tilde{M}}\,,\,\, {\rm where}\nonumber
\end{align}\vspace{-2ex}
\begin{align*}\begin{split}
\tilde{M}& :=\frac{1}{2(1-\tilde{\alpha})}\left(\tilde{s}^2(1-\tilde{\alpha})^2+
|(1-\tilde{\alpha})\tilde{y}+\tilde{u}|^2\right)=\\
&=\frac{\tilde{s}^2+|\tilde{y}|^2+1}{2}-\tilde{\alpha}\frac{\tilde{s}^2+|\tilde{y}|^2-1}{2}+\tilde{u}^t\tilde{y}
\end{split}\end{align*}
\end{lem}
\noindent
{\em Proof:} Direct computation using (\ref{def-B}) and (\ref{def-C}). \\
\dowl

\noindent
By   direct computations one also verifies that  (\ref{eq-iwasawa}) implies equality 
\notka{eq-orbit1}
\begin{equation*}\label{eq-orbit1}
\tilde{s}(\tilde{\alpha}-1) =\frac{\alpha-1}{s}
\end{equation*}
Clearly, it  follows that $\alpha=1\iff\tilde{\alpha}=1$ and  then:
\begin{equation*}
u=w=\tilde{u}=\tilde{w}=0\,,\,\,\Lambda=\tilde{\Lambda}\,,\,\,s=\tilde{s}\,,\,\,\tilde{y}=\Lambda y
\end{equation*}
For  $\alpha\neq1$ (or, equivalently, $\tilde{\alpha}\neq 1$) we have:  \notka{eq-orbit2}
\begin{equation}\label{eq-orbit2}
\Lambda-\frac{1}{\alpha-1}u w^t =\tilde{\Lambda}-\frac{1}{\tilde{\alpha}-1}\tilde{u} \tilde{w}^t
\end{equation}

Moreover, for two given  elements of $B$: $(\Lambda,u,w,\alpha)$ and $(\tilde{\Lambda},\tilde{u},\tilde{w},\tilde{\alpha})$ with $\alpha\neq 1,\tilde{\alpha}\neq 1$ 
that satisfy  (\ref{eq-orbit2}), there exist (not unique) $(s,y)\in C$ and  $(\tilde{s}, \tilde{y})\in C$ such that equation  (\ref{eq-iwasawa}) is fulfilled; they are given by:
\begin{align*}
\tilde{s} s &=\frac{\alpha-1}{\tilde{\alpha}-1}\,,\,\,&\,
\tilde{y} &=\frac{u}{s(1-\tilde{\alpha})}-\frac{\tilde{u}}{1-\tilde{\alpha}}\,,&\,y &=\frac{\tilde{w}}{\tilde{\alpha}-1}-
\frac{s w}{\alpha-1}
\end{align*}

The  decomposition (\ref{eq-iwasawa}) defines  the groupoid $G_B: G\rightrightarrows B$.
We  will write $(\Lambda,u,w,\alpha;s,y)$ for the product $(\Lambda,u,w,\alpha)(s,y)$ and use $(\Lambda,u,w,\alpha;s,y)$ to denote elements of $G_B$.

The multiplication relation is given by:\notka{GB-multip}
\begin{equation*}\label{GB-multip}
\begin{split}
m_B:=\left\{\left(\Lambda,u,w,\alpha;s s_1,s_1 y+y_1;;\,
\Lambda,u,w,\alpha;s,y\,;;\,\tilde{\Lambda},\tilde{u},\tilde{w},\tilde{\alpha};s_1,y_1\right): \right.\\
\left.(\Lambda,u,w,\alpha;s,y)=(\tilde{s},\tilde{y})(\tilde{\Lambda}, \tilde{u}, \tilde{w}, \tilde{\alpha})\right\}\subset G_B\times G_B\times G_B, 
\end{split}
\end{equation*}
%%%%%%%%%%%%%%
%%%%%%%%%%%%%%%%%%%%%%%%%%%%%%%%%%%%%
We will use  a detailed structure of this groupoid to verify (3) of Assumptions \ref{basic-assumpt}; the structure is summarized in the following lemma.
\begin{lem}\label{lemma-GB-struct}\notka{lemma-GB-struct}
\begin{enumerate}
\item The isotropy group of  $(\Lambda,0,0,1)\in B$ 
(i.e. of elements of $SO(n)$ embedded in $SO(n+1)$ in the upper left corner) is  equal  to $\{(\Lambda,0,0,1;s,y): (s,y)\in C)\}\simeq C$ 
and for $(\Lambda,u,w,\alpha)\in B\,,\,\alpha\neq 1$ is one dimensional $\{(\Lambda,u,w,\alpha;s, \frac{s-1}{1-\alpha} w ):s\in \R_+\}\simeq \R_+$
\item $G_B$ is a disjoint  union of an open groupoid $\Gamma_0$ over $SO(n+1)\setminus SO(n)$ and 
a group bundle $\Gamma_1:=\{(\Lambda,0,0,1;s,y): \Lambda\in SO(n),\,(s,y)\in C \} \simeq  SO(n)\times C$; %($SO(n)$ is embedded into $SO(n+1)$ into left upper corner);
\item The open groupoid $\Gamma_0$ is a product $\Gamma_0=O(n)^-\times \Gamma_{00}$ of a manifold-groupoid $O(n)^-$  and a transitive groupoid $\Gamma_{00}$, where
$O(n)^-$ stands for the component of the orthogonal group with a  negative determinant;
\item The transitive groupoid $\Gamma_{00}$ is a product $\Gamma_{00}=\R_+\times(\R^n\times\R^n)$ of a group and a pair groupoid.
\end{enumerate}
\end{lem}
\noindent
{\em Proof:} 1) and 2) are direct consequences of formulae  (\ref{solution-iwasawa1}) and  (\ref{solution-iwasawa2});\\
3) For $v\in\R^n$ let $R(v)$ be the orthogonal reflection in $\R^{n+1}$  along the direction of $(v,-1)^t$ i.e. 
$$R(v)=\left(\begin{array}{cc} I-\frac{2 v v^t}{1+|v|^2} & \frac{2 v}{1+|v|^2} \\
\frac{2 v^t}{1+|v|^2} &\frac{|v|^2-1 }{1+|v|^2}\end{array}\right)$$
%%%%%%%%%%%%%%%
Let us consider the map: 
%%%%%%%%%%%%%%%%
$$\Phi: O(n)^-\times \R^n\ni (K,v)\mapsto \left(\begin{array}{cc} K &0\\0&1\end{array}\right) R(v) 
\in SO(n+1)$$
Clearly  $\Phi(K,v) \in SO(n+1)\setminus SO(n)$. 
Moreover, using (\ref{param-SOn}),  it is easy to see that for $\alpha\neq 1$ the matrix 
$\Lambda-\frac{1}{\alpha-1}u w^t$ is orthogonal and
$$\left(\begin{array}{cc} \Lambda& u\\w^t &\alpha\end{array}\right) 
\left(\begin{array}{cc} I& 0\\\frac{w^t}{1-\alpha}  & 1\end{array}\right)
\left(\begin{array}{cc} I& -w\\0 & 1\end{array}\right)=
\left(\begin{array}{cc} \Lambda-\frac{1}{\alpha-1}u w^t& 0\\ \frac{w^t}{1-\alpha}&-1\end{array}\right),$$
therefore $\det(\Lambda-\frac{1}{\alpha-1}u w^t)=-1$ and we can define 
$$\Psi:SO(n+1)\setminus SO(n)\ni \left(\begin{array}{cc} \Lambda& u\\w^t &\alpha\end{array}\right) \mapsto 
\left( \Lambda-\frac{1}{\alpha-1}u w^t,\frac{w}{1-\alpha}\right)\in  O(n)^-\times \R^n.$$
By direct computation one verifies that $\Psi=\Phi^{-1}$. Clearly,  both mappings are smooth, so both are diffeomorphisms.
%%%
By (\ref{solution-iwasawa1}) we obtain:
$$\Phi(K,v)\,(s,y)=(\tilde{s},\tilde{y})\,\Phi(K, sv-y)$$
$$\tilde{s}=\frac{1+|s v-y|^2}{s(1+|v|^2)}\,,\,\,\,\,
\tilde{y}=K\left(y-\frac{2v}{1+|v|^2}(v^t y +\frac{s^2-|y|^2-1}{2 s})\right)$$
These formulae  show that $\Gamma_{00}:=\R^n\times C$ is a (right) transformation groupoid with the action:
$\R^n\times C\ni (v ;s,y)\mapsto s v-y\in \R^n$. It is clear that this action is transitive.\\
4)  We will use the following general but simple:
\begin{lem}\label{lemma-transitive} \notka{lemma-transitive}Let $\Gamma\rightrightarrows E$ be a transitive groupoid. For $e_0\in E$ 
let $p: E\ni e \mapsto p(e)\in e_L^{-1}(e_0)$ be a section of the right projection, such that $p(e_0)=e_0$.
Let $G$ be an isotropy group of $e_0$. 
The map $G\times E\times E\ni(g,e_1,e_2)\mapsto s(p(e_1)) g p(e_2)\in\Gamma$ is an isomorphism of groupoids.
($G\times E\times E$ is a product of a group and pair groupoid).
\end{lem}
The application of  the lemma %(\ref{lemma-transitive}) 
(choose $v_0=0$) gives us a groupoid isomorphism:
$$\R_+\times\R^n\times \R^n\ni(s;x_1,x_2)\mapsto (x_1; s, s x_1-x_2)\in \R^n\times C,$$
which in our situation clearly is a diffeomorphism.\\
\dowl

To find the set $B'$ i.e. $B\cap CA$ we have to solve the equation:\notka{eq-B-prim}
\begin{equation}\label{eq-B-prim}
(z,U,d)=(s,y)(\Lambda,u,w,\alpha)\,,\,\,(z,U,d)\in A\,,\,(s,y)\in C\,,\,(\Lambda,u,w,\alpha)\in B
\end{equation}
Using (\ref{def-zUd}, \ref{def-B}) and (\ref{def-C}), by direct computation, one verifies:

\begin{lem}\label{lemma-B-prim}\notka{lemma-B-prim}
(a) For  $(z,U,d)\in A$  solutions $(s,y)\in C\,,\,(\Lambda,u,w,\alpha)\in B$ of (\ref{eq-B-prim}) are given by:
\begin{align}\label{solution-B-prim-1}
s& =\frac{1+|z|^2}{1-|z|^2}\quad, & y&=\frac{- 2 z}{1-|z|^2}\quad,\nonumber\\
\Lambda&= (I-\frac{2 z z^t}{1+|z|^2})U D\quad, &
u&=\frac{- 2 d z }{1+|z|^2}\quad,& w&=\frac{2 D U^t z}{1+|z|^2}\quad, & \alpha&=\frac{d(1-|z|^2)}{1+|z|^2}. \notka{solution-B-prim-1}
\end{align}
(b) For $(\Lambda,u,w,\alpha)\in B$ with $\alpha\neq 0$, solutions $(s,y)\in C\,,\,(z,U,d)\in A$ of (\ref{eq-B-prim}) are given by:
\begin{align}\nonumber\label{solution-B-prim-2}
s&=\frac{1}{|\alpha|}\quad,& y& =\frac{u}{\alpha}\quad,\\
z& =-\frac{\sign (\alpha) u}{1+|\alpha|}\quad,& U&=(\Lambda-\frac{\sign(\alpha)}{1+|\alpha|} u w^t)D\quad, & d& =\sign(\alpha). \notka{solution-B-prim-2}
\end{align}
\end{lem}
\dowl

\noindent
This way  we obtain $B'$ and projections $\tilde{c}_L$ and $a_R$ defined by decomposition $CA\ni g=\tilde{c}_L(g)a_R(g)$ (restricted to $B'$).
\begin{align}
B'& =\{(\Lambda,u,w,\alpha)\in B: \alpha\neq 0\}\label{def-B-prim}\notka{def-B-prim}\\
\tilde{c}_L(\Lambda,u,w,\alpha)& =(|\alpha|,-sgn(\alpha) u )\label{def-tcL}\notka{def-tcL}\\
a_R(\Lambda,u,w,\alpha)& =\left(-\frac{\sign (\alpha) u}{1+|\alpha|}, \,(\Lambda-\frac{\sign(\alpha)}{1+|\alpha|} u w^t)D, \,\sign(\alpha)\right)\label{def-aR} \notka{def-aR}
\end{align}
Finally, using the definition (\ref{def-twist}) and  the formula above we obtain the twist $T$:
\begin{equation*}\label{def-twist-kappa}
T=\{(\Lambda_1,u_1,w_1,\alpha_1;\frac{1}{|\alpha_2|},\frac{u_2}{\alpha_2};;\Lambda_2,u_2,w_2,\alpha_2;1,0)\in G_B\times G_B\,:\,\alpha_2\neq 0\}\notka{def-twist-kappa}
\end{equation*}
%%%%%%%%%%%%%%%%%%%%%%%%%%%%%%%%%%%%%%%%%%%%%%%%%%%%%%%%%%%%%%%%%%%%%%%%%%%%%%%%%%%%%%%%%%%%%%%
%%%%%%%%%%%%%%%%%%%%%%%%%%%%%%%%%%%%%%%%%%%%%%%%%%%%%%%%%%%%%%%%%%%%%%%%%%%%%%%%%%%%%%%%%%%%%%
\section{The $C^*$-algebra}

In this section we will use Prop.~\ref{prop-Delta} to get the $C^*$-algebra and the comultiplication.  
To this end, conditions listed in assumptions \ref{basic-assumpt} have to be verified. 
The first one is clear and the second was proven in  \cite{PS-triple}, it remains to check $(3)$ and $(4)$.

\noindent
{\bf\small Condition \ref{basic-assumpt} (3) \,} Let $U:=b_L^{-1}(B')$ and $\sA(U)$ be the linear space of elements from $\sA(G_B)$ supported in $U$.
We will show that $\sA(U)$ is dense in $C^*_r(G_B)$.\\
Let $B'$ be as in (\ref{def-B-prim}) and define  $B'':=\{(\Lambda,u,w,\alpha)\in B: \alpha\neq 1\}$, then $B=B'\cup B''$ i.e.  $\{B',B''\}$ is an open cover of $B$ and $b_L^{-1}(B'')=\Gamma_0$.
Let $\tilde{\chi}',\tilde{\chi}''$ be a partition of unity subordinated to the cover $\{B',B''\}$, 
$\chi'(\gamma):=\tilde{\chi}'(b_L(\gamma))$ and  $\chi''(\gamma):=\tilde{\chi}''(b_L(\gamma))$.
For $\om\in\sA(G_B): \om=\chi'\om+\chi''\om$ and $\chi'\om\in\sA(U)\,,\,\chi''\om\in\sA(\Gamma_0)$; so to prove that 
$\sA(U)$ is dense w $C^*_r(G_B)$ is sufficient to show that any element in $\sA(\Gamma_0)$ can be approximated by elements from $\sA(U)$, 
in particular from $\sA(U\cap\Gamma_0)$.
In this way we can transfer the whole problem  to $\Gamma_0$ and use its structure described in lemma \ref{lemma-GB-struct};
 in this presentation $U\cap\Gamma_0=\{(K,t, x_1,x_2)\in O(n)^-\times \R_+\times \R^n\times \R^n\, :|x_1|\neq 1\}$.
 Let $dK$ denotes a measure on $O(n)^-$ defined by (the square of) some smooth, positive, non vanishing half-density; $dx_1$ be 
the Lebesgue measure on $\R^n$ and $\frac{d s}{s}$ be the Haar measure on $\R_+$. Then $L^2(\Gamma_0)$ can be identified with 
$L^2\left(O(n)^-\times \R_+\times \R^n\times \R^n; dK \frac{d s}{s} d x_1 d x_2\right)$ and $\sA(\Gamma_0)$ acts on $L^2(\Gamma_0)$ by:\notka{def-rep-id-GB}
\begin{equation}\label{def-rep-id-GB}
(\pi(f)\Psi)(K,t,x_1,x_2):=(f *\Psi)(K,t,x_1,x_2):=\int \frac{d s}{s}\, d y\, f(K,s,x_1,y)\Psi(K,t/s,y,x_2);
\end{equation}
note that it is sufficient to consider $\Psi\in\sD(\Gamma_0)$.
We will prove the following estimate:
\notka{lemma-basic-estimate}
\begin{lem}\label{lemma-basic-estimate} Let $f$ be a continuous function supported in a product of compact sets 
$L\times M\times N\times R \subset O(n)^-\times \R_+\times \R^n\times \R^n$. 
Then the norm of the operator $\pi(f)$ given by (\ref{def-rep-id-GB}) satisfies:
 \begin{equation}\label{eq-basic-estimate}||\pi(f)||\leq \sup|f|\, \nu(M)\,\sqrt{\mu(N)\mu(R)},
\end{equation}
where $\nu$ denotes the Haar measure on $\R_+:\nu(M):=\int_M \frac{d s}{s}$ and $\mu$ the Lebesgue measure on $\R^n$.
\end{lem}

For $\epsilon > 0$ let  $\sO_\epsilon$ be  an open neighborhood of the unit sphere in $\R^n$ with $\mu(\sO_\epsilon)\leq\epsilon$; 
let $\tilde{\chi}_\epsilon$ be a smooth function supported in $\sO_\epsilon$ such that 
$0 \leq \tilde{\chi}_\epsilon \leq 1$, $\tilde{\chi}_\epsilon =1$ on the unit sphere and 
the function $\chi_\epsilon$ on $\Gamma_0$ be defined by $\chi_\epsilon(K,t,x_1,x_2):=\tilde{\chi}_\epsilon (x_1)$.

 Let $f\in\sD(\Gamma_0)$ be supported in $L\times M\times N\times R$.  
We have  $f=(f-\chi_\epsilon f)+\chi_\epsilon f$ and $(f-\chi_\epsilon f )\in \sD(U\cap\Gamma_0)$. 
By the lemma \ref{lemma-basic-estimate} 
$$||\pi(\chi_\epsilon f)||\leq \sup|f| \nu(M) \sqrt{\mu(R)} \sqrt{\mu(\sO_\epsilon\cap N)}\leq \sup|f| \nu(M) \sqrt{\mu(R)} \sqrt{\epsilon} $$
So really $f$ can be approximated by elements from $\sD(U\cap\Gamma_0)$. It remains to prove lemma \ref{lemma-basic-estimate}.

{\em Proof of lemma \ref{lemma-basic-estimate}:} 
We  just apply the Schwartz inequality several times.
Let $\Psi$ be smooth and compactly supported; by (\ref{def-rep-id-GB}): 
$$(\Psi|f*\Psi)=\int dK \frac{d t}{t} d x_1 d x_2\, \overline{\Psi(K,t,x_1,x_2)}\int\frac{d s}{s} d y\,f(K,s,x_1,y)\Psi(K,t/s,y,x_2)$$
and 
\begin{equation}\label{norm-est}
|(\Psi|f*\Psi)|\leq \int dK \frac{d t}{t} d x_1 d x_2\,|\Psi|(K,t,x_1,x_2)\int\frac{d s}{s} d y\,|f|(K,s,x_1,y)|\Psi|(K,t/s,y,x_2)
\end{equation}

We write the integral as iterated integral: $\int dK d x_2 \int d x_1 \int d y \,\int \frac{d t}{t}\,\int\frac{d s}{s}$ 

For fixed $(K,y,x_1,x_2)$ let us  define functions: $\displaystyle\Psi_1:\R_+\ni t\mapsto \Psi_1(t):=|\Psi(K,t,x_1,x_2)|\,$
$$f_1:\R_+\ni t\mapsto f_1(t):=|f(K,t,x_1,y)|\,\,;\,\,\Psi_2:\R_+\ni t\mapsto \Psi_2(t):=|\Psi(K,t,y,x_2)|$$

With these definitions we have the estimate: 
$$\int \frac{d t}{t}\,|\Psi|(K,t,x_1,x_2)\,\int\frac{d s}{s}|f|(K,s,x_1,y)|\Psi|(K,t/s,y,x_2)=|(\Psi_1|f_1 *\Psi_2)|
\leq ||\Psi_1||_2||\Psi_2||_2||f_1||_1,$$
where the scalar product is in  $L^2(\R_+,\frac{d s}{s})$ and  $||\cdot||_2$ norms refer to this space, $*$ is the convolution in $\R_+$ and  $||f||_1$ is $L^1$ norm; 
these norms are continuous functions of remaining variables.\\
Let us now define (continuous, compactly supported) functions 
$ \tilde{\Psi}_1,\,\tilde{\Psi}_2,\, \tilde{f}_1 : O(n)^-\times \R^n\times \R^n\rightarrow \R $: 
$$\tilde{\Psi}_1(K,x_1,x_2):=||\Psi_1||_2=\left[\int\,\frac{d t}{t} |\Psi(K,t,x_1,x_2)|^2\right]^{1/2},\,\tilde{\Psi}_2(K,x_2,y):=||\Psi_2||_2$$  and  
$$\tilde{f}_1(K,x_1,y):=||f_1||_1=\int t\,\frac{d t}{t} |f(K,t,x_1,y)|.$$
 The right hand side of (\ref{norm-est}) is estimated by:
$$\int dK d x_2\int d x_1\,\tilde{\Psi}_1(K,x_1,x_2)\,\int d y \,\tilde{f}_1(K,x_1,y)\tilde{\Psi}_2(K,x_2,y)$$
By the Schwartz inequality for $y$-integration we  get an estimate for the integral above by:
\begin{equation}\label{norm-est1}
\int dK d x_2 \int d x_1 \tilde{\Psi}_1(K,x_1,x_2) \left[\int d y \,(\tilde{f}_1(K,x_1,y))^2\right]^{1/2}\,
\left[\int d y \,(\tilde{\Psi}_2(K,x_2,y))^2\right]^{1/2}
\end{equation}

Again let us denote $\tilde{f}_2(K,x_1):=\left[\int d y \,(\tilde{f}_1(K,x_1,y))^2\right]^{1/2}$ and 
$\tilde{\Psi}_3(K,x_2):=\left[\int d y \,(\tilde{\Psi}_2(K,x_2,y))^2\right]^{1/2}$

The integral (\ref{norm-est1}) reads
$$\int dK d x_2 \, \tilde{\Psi}_3(K,x_2)\,\int d x_1 \, \tilde{f}_2(K,x_1)\,\tilde{\Psi}_1(K,x_1,x_2).$$

Schwartz inequality again and we can estimate it  by:
$$\int dK d x_2 \, \tilde{\Psi}_3(K,x_2)\,\tilde{\Psi}_4(K,x_2)\tilde{f}_3(K),$$
where $\tilde{\Psi}_4(K,x_2):=\left[\int d x_1\,(\tilde{\Psi}_1(K, x_1,x_2))^2\right]^{1/2}$ and   
$\tilde{f}_3(K):=\left[\int d x_1\,(\tilde{f}_2(K, x_1))^2\right]^{1/2}$.\\
And, finally, this  integral we estimate by 
$$\sup|\tilde{f}_3|\,\left[\int dK d x_2 (\tilde{\Psi}_4(K,x_2))^2\right]^{1/2}\,
\left[\int dK d x_2 (\tilde{\Psi}_3(K,x_2))^2\right]^{1/2}$$
(again the Schwartz inequality was used).

But 
$$\int dK d x_2 (\tilde{\Psi}_4(K,x_2))^2=\int dK d x_2 \int d x_1\,(\tilde{\Psi}_1(K, x_1,x_2))^2=$$
$$=\int dK d x_2 \int d x_1\,\int \frac{d t}{t}\,|\Psi|^2(K,t,x_1,x_2)=\|\Psi\|^2$$
and, in a similar way,
$$\int dK d x_2 (\tilde{\Psi}_3(K,x_2))^2=
\|\Psi\|^2$$
Thus we get an inequality $|(\Psi|f*\Psi)|\leq \sup|\tilde{f}_3| \|\Psi\|^2$, therefore  $||\pi(f)||\leq \sup|\tilde{f}_3|$.

$$|\tilde{f}_3(K)|^2=\int d x_1 (\tilde{f}_2(K,x_1))^2=\int d x_1\int d y\left[\int \frac{d s}{s} |f(K,s,x_1,y)|\right]^2\leq$$
$$\leq\int_N d x_1\int_R d y \,(\sup|f|)^2\nu^2(M)=(\sup|f|)^2\nu^2(M) \mu(N) \mu(R)$$
This is the estimate (\ref{eq-basic-estimate}) and the lemma is proven.
\\\dowl

\noindent
This way by Prop. \ref{prop-Delta} (a) we obtain the equality  of $C^*$-algebras: 
\begin{equation}\cred(\Gamma_A)=\cred(\G_{B'})=\cred(G_B).\end{equation}

\noindent
To get the twist and the comultiplication we have to verify:\\
{\bf\small Condition \ref{basic-assumpt} (4) \,} %{\color{red} te rachunki sprawdzone }
For $M>1$ and $1>\delta>\epsilon>0$ let us define  compact sets $K_M\subset C, K_\delta\subset B'$ and 
 an open neighborhood $V_\epsilon$ of $B\setminus B'$ in $B$:
$$K_M:=\{(s,y)\in C: \frac{1}{M}\leq s\leq M\,,\,|y|\leq M\}\,,\,\,K_\delta:=\{(\Lambda,u,w,\alpha)\in B: |\alpha|\geq \delta\}\,,$$
$$V_\epsilon:=\{(\Lambda,u,w,\alpha)\in B : |\alpha| <  \epsilon\}.$$
We will show that for any $M>1$ any $0<\delta<1$ and any $0< \epsilon< \delta$ there exists $\epsilon'$ such that $\mu(B,K_\delta,K_M;V_{\epsilon'})<\epsilon$ (notation as in 
Assumptions \ref{basic-assumpt}.)  Since any compact in $B'$ is contained in some $K_\delta$ and any compact in  $C$ is contained in some $K_M$, this is sufficient.

Let  $b=(\Lambda, u, w, \alpha)$ and $b_1=(\Lambda_1,u_1,w_1,\alpha_1)\in K_\delta$; we want to find the set  
$\displaystyle Z(b,b_1,K_M;V_\epsilon):=K_M\cap \{c\in C: b_R(b c)b_1\in V_\epsilon\}.$
Let  $c=(s,y)\,,\,\,b_R(b c)=:(\tilde{\Lambda},\tilde{u},\tilde{w},\tilde{\alpha})$ and $b_R(b c) b_1=:(\Lambda_2,u_2,w_2,\alpha_2)$. 

\noindent Then $\alpha_2=\tilde{w}^t u_1 +\tilde{\alpha}\alpha_1$ and 
using solutions of eq. (\ref{eq-iwasawa}) we get:\\ 

$\displaystyle \tilde{\alpha}=\left\{\begin{array}{cc}
    \frac{|r|^2-1}{|r|^2+1} & \alpha\neq 1\\1 &
    \alpha=1\end{array}\right.\,$, 
$\,\displaystyle\tilde{w}=\left\{\begin{array}{cc} \frac{2 r}{1+|r|^2} &
    \alpha\neq 1\\ 0& \alpha=1\end{array}\right.\,$, 
$\,\displaystyle\alpha_2=\left\{\begin{array}{cc}\frac{2 r^t
      u_1+\alpha_1(|r|^2-1)}{1+|r|^2} & \alpha\neq 1\\ \alpha_1 &
    \alpha=1\end{array}\right.$,\\

\noindent  where $r:=\frac{s w}{1-\alpha} -y$.

\noindent
Now we  solve for $(s,y)$ the inequality   $|\alpha_2|< \epsilon$ with the additional assumption   $0<\epsilon<\delta$; 

\noindent
There is no solution for  $\alpha=1$ and for $\alpha\neq 1$ we have:

$$-\epsilon(1+|r|^2)<2 r^t u_1+\alpha_1(|r|^2-1)<\epsilon(1+|r|^2)$$
after some manipulation we get for $\alpha_1>0$:
$$|r+\frac{u_1}{\alpha_1-\epsilon}|^2< \frac{1-\epsilon^2}{(\alpha_1-\epsilon)^2}\,\,\,{\rm and}\,\,\,
|r+\frac{u_1}{\alpha_1+\epsilon}|^2> \frac{1-\epsilon^2}{(\alpha_1+\epsilon)^2}$$
and for $\alpha_1<0$:
$$|r+\frac{u_1}{\alpha_1-\epsilon}|^2> \frac{1-\epsilon^2}{(\alpha_1-\epsilon)^2}\,\,\,{\rm and}\,\,\,
|r+\frac{u_1}{\alpha_1+\epsilon}|^2< \frac{1-\epsilon^2}{(\alpha_1+\epsilon)^2}.$$

\noindent
Both situations can be described uniformly as: 
 
$$\left|r+\frac{sgn(\alpha_1)}{|\alpha_1|-\epsilon} u_1 \right|^2<\frac{1-\epsilon^2}{(|\alpha_1|-\epsilon)^2}\,\,\,{\rm
  and}\,\,\,
\left|r+\frac{sgn(\alpha_1)}{|\alpha_1|+\epsilon} u_1\right|^2>
\frac{1-\epsilon^2}{(|\alpha_1|+\epsilon)^2}$$
or in terms of $(s,y)$:

$$\left|\frac{s w}{1-\alpha}
  +\frac{sgn(\alpha_1)}{|\alpha_1|-\epsilon} u_1 - y \right|^2<\frac{1-\epsilon^2}{(|\alpha_1|-\epsilon)^2}\,\,\,{\rm
  and}\,\,\,
\left|\frac{s w}{1-\alpha}+\frac{sgn(\alpha_1)}{|\alpha_1|+\epsilon}
  u_1 - y \right|^2>
\frac{1-\epsilon^2}{(|\alpha_1|+\epsilon)^2}$$
For fixed $s$ this is the intersection of the (larger) ball
centered at $y_1:=\frac{s w}{1-\alpha}
  +\frac{sgn(\alpha_1)}{|\alpha_1|-\epsilon} u_1$ with a radius
  $r_1:=\frac{\sqrt{1-\epsilon^2}}{|\alpha_1|-\epsilon}$ with the
  exterior of the (smaller) ball centered at $y_2:=\frac{s w}{1-\alpha}
  +\frac{sgn(\alpha_1)}{|\alpha_1|+\epsilon} u_1$ with a radius
  $r_2:=\frac{\sqrt{1-\epsilon^2}}{|\alpha_1|+\epsilon}$.  Because of the inequality 
$$ |y_1-y_2|=\left(\frac{1}{|\alpha_1|-\epsilon}- \frac{1}{|\alpha_1|+\epsilon}\right) \sqrt{1-\alpha_1^2}
< \left(\frac{1}{|\alpha_1|-\epsilon}- \frac{1}{|\alpha_1|+\epsilon}\right)\sqrt{1-\epsilon^2}=r_1-r_2,$$
the smaller ball is contained in the larger one, and 
the volume of this intersection is equal to:
$$F(n)(r_1^n-r_2^n)=F(n)\frac{(1-\epsilon^2)^{n/2}}{(|\alpha_1|-\epsilon)^n}
\left(1-\left(1-\frac{2 \epsilon}{|\alpha_1|+\epsilon}\right)^n\right)
\leq F(n) \frac{(1-\epsilon^2)^{n/2}}{(|\alpha_1|-\epsilon)^n} \frac{2 n\epsilon}{|\alpha_1|+\epsilon}\leq$$
$$ \leq F(n) \frac{1}{(\delta-\epsilon)^n}\frac{2 n\epsilon}{\delta}\leq 2 n F(n)\frac{\epsilon}{\delta (\delta-\epsilon)^n },$$
where $F(n)r^n$ is a volume of n-dimensional ball. 
%what converges uniformly (for $|\alpha_1|>\delta>\epsilon$) to $0$.
In this way we obtain
$$\mu(b,b_1,K_M;V_\epsilon)=\int_{Z(b,b_1,K_M;V_\epsilon)}\frac{ds}{s} dy\leq \epsilon \,\frac{\log M}{\delta} \frac{4 n F(n) }{(\delta-\epsilon)^n}$$
i.e. $$\mu(B,K_\delta,K_M;V_\epsilon)\leq \epsilon \,\frac{\log M}{\delta} \frac{4 n F(n) }{(\delta-\epsilon)^n}$$
The right hand side  goes to $0$ as $\epsilon\rightarrow 0$, and the fourth condition of (\ref{basic-assumpt}) is satisfied.\\\dowl

Now, by the Prop.~\ref{prop-Delta} (b), (c),  we get the comultiplication $\Delta$ on $\cred(\G_{B'})=\cred(G_B)$ satisfying the density condition (\ref{density-c}).

%%%%%%%%%%%%%%%%%%%%%%%%%%%%%%%%%%%%%%%%%%%%%%%%%%%%%%%%%%%%%%%%%%%%%%%%%%%%%%%%%%%%%%%%%%%%%
%%%%%%%%%%%%%%%%%%%%%%%%%%%%%%%%%%%%%%%%%%%%%%%%%%%%%%%%%%%%%%%%%%%%%%%%%%%%%%%%%%%%%%%%%%%%%
\section{Generators and relations}
%%%%%%%%%%%%%%%%%%%%%%%%%%%%%%%%%%%%%%%%%%%%%%%%%%%%%%%%%%%%%%%%%%%%%%%%%%%%%%%%%%%%%
In the previous section the $C^*$-algebra of the quantum $\kappa$-Poincar\'{e} Group together with comultiplication was defined. In this section we describe its generators, 
commutation relations among them and look closer at the twist.

\subsection{General formulae}
Let  $ u $ be  a bisection of a differential groupoid $\G\rightrightarrows E$ and $f$ a smooth, bounded function on $E$. They define bounded
operators on $L^2(\Gamma)$, denoted by $\hat{u}$ and $\hat{f}$, which are multipliers of $C^*_r(\Gamma)$: $\hat{u}$ acts by a push-forward of half-densities and the action of 
$\hat{f}$ is defined by 
$(\hat{f} \psi)(\gamma):=f(\el(\gamma)) \psi(\gamma)\,,\,\psi\in\omh(\G)$ \cite{DG}. These operators satisfy: 
\notka{bis-fun-komut1}
\begin{equation}\label{bis-fun-komut1}
\hat{u} \hat{f}= \widehat{f_{\scriptsize u}} \hat{u}\,\,,\,\,{\rm where\,\,} f_{\scriptsize u}(e):=f(e_L(u^{-1}e))\,,\quad e\in E  
\end{equation}
In particular, for a one-parameter group of bisections $u_t$,  we have \notka{bis-fun-komut2}
\begin{equation}\label{bis-fun-komut2}
\widehat{u_t} \hat{f} \widehat{u_{-t}} = \widehat{f_t}\,\,,\,\,{\rm where }\,\, f_t(e):=f(e_L(u_{-t} e)) 
\end{equation}
%%%%%%%%%%%%%%%%%%%%%%%%%%
For a DLG  $(G;B, C)$ and $c_0\in C$,  by  $B c_0$ we denote  the bisection of $G_B$ defined as $B c_0:=\{b c_0 :b\in B\}$. It acts   on $G_B$  by 
\notka{Bc0-action0}
\begin{equation} \label{Bc0-action0} (B c_0)(b  c )= b_R(b c_0^{-1}) c_0 c = c_L^{-1}(b  c_0^{-1}) b  c= c_R(c_0 b^{-1}) b c  \end{equation}
%%%%%%%%%%%%%%%%%%%%%%%%%%%%%%%%%%%%%%%%%%%%%%%%%%%%%%%%%%%%%%%%%%%%%%%%%%%%%%%%%%%%%%%%%%%%%%%%
\noindent
$G_B$ is a (right) transformation groupoid $B\rtimes C$ for  the action $B\times C \ni (b,c)\mapsto b_R(bc)\in B$. The $C^*$-algebra $C^*_r(G_{B})$ is the reduced crossed product
$C_0(B)\rtimes_r C$ (see e.g \cite{PSAx}, prop.~5.2, where this identification is described in details). If $C$ is {\em an amenable group}, and this is our case 
(or more generally, this is the case of $C=AN$ coming from the Iwasawa decomposition $G=K(AN)$) 
{\em reduced and universal  crossed products coincide} (\cite{DW}, Thm 7.13, p. 199). 
Moreover, since the universal $C^*$-algebra of a differential groupoid $C^*(G_{B})$ as defined in \cite{DG} is, 
for transformation groupoids,  ``something between'' universal and reduced crossed products, for amenable groups $C$  we have $C^*(G_{B})=C^*_r(G_{B})$; thus
any morphism of  differential groupoids $h:G_{B} \rel \Gamma$ defines a $C^*$-morphism of corresponding reduced and  universal algebras.
\begin{re}
Thus our $C^*(G_B)$ has $SO(n)$ family of ``classical points'' i.e. characters given by $0$-dimensional orbits of $G_B$ described in lemma \ref{lemma-GB-struct}.
\end{re}

The canonical morphisms $i_{C}\in Mor(C^*(C), C^*(G_{B}))$ and $i_{B}\in Mor(C_0(B),C^*(G_{B}))$ are given by (extensions of) the following actions on $\sA(G_{B})$:
$$ i_{C}(c_0)\omega:= (B c_0) \omega\,,\quad (i_{B}(f)\omega)(g):= f(b_L(g)) \omega(g)\,,\quad\quad c_0\in C, g\in G,  f\in C^\infty_0(B),\omega\in \sA(G_{B})$$
%%%%%%%%%%%%%%%%%%%%%%%%%
If $X_1,\dots, X_m$ are generators  of $C^*(C)$,  $f_1,\dots, f_k$ are generators of $C_0(B)$ (in the sense of \cite{Wor2})  and $C$ is amenable then 
(by the universality of crossed product) 
$C^*(G_{B})$ is generated by  $i_{C}(X_1),\dots$, $i_{C}(X_m)$, $i_{B}(f_1),\dots$, $i_{B}(f_k)$.

In our situation, it is clear that $C(B)$ is generated by matrix elements $(\Lambda,u,w,\alpha)$ in (\ref{def-B}) and, 
by the results of \cite{Wor2},  $C^*(C)$ is generated by any basis in $\gotc$. Thus we have
%%%%%%%%%%%%%%%%
\begin{prop} Let groups $B$ and $C$ be defined by (\ref{def-B}, \ref{def-C}); let  $(X_0,\dots,X_n)$ be a basis in $\gotc$ and $(\Lambda,u,w,\alpha)$ be matrix elements of $B$. 
Elements $(i_C(X_0), \dots,  i_C(X_n), i_B(\Lambda), i_B(u), i_B(w), i_B(\alpha))$ are generators of $C^*(G_B)$.\\
\dowl
\end{prop}

%%%%%%%%%%%%%%%%%%%%
For $\kropka{c}\in\gotc$,  let $u_{\kropka{c}}$ and $X^r_{\kropka{c}}$ be, respectively,  the one-parameter group of bisections and the  right invariant vector field on $G_{B}$  defined by
\notka{def-u-Xr-ckropka}
\begin{equation}\label{def-u-Xr-ckropka}  u_{\kropka{c}}(t):=B\exp(t \kropka{c}) \quad,\quad X^r_{\kropka{c}}(b c ):=\left.\frac{d}{d t}\right|_{t=0} u_{\kropka{c}}(t)(b c)=( \adc(b)(\kropka{c})) b c , 
\end{equation}
where $\adc(g)$ is defined in (\ref{adb-adc}) (the last equality follows easily from (\ref{Bc0-action0})).

Let $\kot (X^r_{\kropka{c}})$ be the projection of $X^r_{\kropka{c}}$ onto  $B$ by $b_L$ i.e.  \notka{def-kotwica}
\begin{equation} \label{def-kotwica}
\kot (X^r_{\kropka{c}})(b):=\left.\frac{d}{d t}\right|_{t=0} b_L(u_{\kropka{c}}(t)(b)) 
\end{equation}
by the straightforward computation, we get:   \notka{kotwica1}
\begin{equation} \label{kotwica1}
\kot(X^r_{\kropka{c}})(b):=- (\rzutB Ad(b) \kropka{c})b
\end{equation}
\begin{prop}\label{prop-generators}\notka{prop-generators} Let $(G;B,C)$ be a DLG and  $\kropka{c}\in\gotc$.  
Let $u_{\kropka{c}}, X^r_{\kropka{c}}$ and $ \kot (X^r_{\kropka{c}})$ be objects defined  in (\ref{def-u-Xr-ckropka}) and (\ref{def-kotwica}).
\begin{enumerate}
\item The one-parameter group $\widehat{u_{\kropka{c}}}$ is strongly continuous and strictly continuous (as a group of multipliers of $C^*_r(G_{B})$). 
Let $A_{\kropka{c}}$ be its generator. The action of $A_{\kropka{c}}$ on $\omh(G_{B})$  and on $\sA(G_{B})$ is given by  $A_{\kropka{c}}=\iota \sL_{X^r_{\kropka{c}}}$ (the Lie derivative).
The linear space $\omh(G_{B})$ is  an essential domain for $A_{\kropka{c}}$ and $A_{\kropka{c}}$ is affiliated to $C^*_r(G_{B})$.
\item For  $\kropka{c},\kropka{e}\in\gotc$, generators $A_{\kropka{c}}, A_{\kropka{e}}$, as  operators  on $\omh(G_{B})$  (or $\sA(G_{B})$),   satisfy:
\begin{equation}\label{komut-A}\notka{komut-A}
[A_{\kropka{c}}, A_{\kropka{e}}]=-\iota A_{[\kropka{c},\kropka{e}]} 
\end{equation}
\item Let $\kropka{c}\in\gotc$ and $f$ be a smooth function on $B$. As operators on $\omh(G_{B})$  (or $\sA(G_{B})$), $A_{\kropka{c}}$ and $\hat{f}$  satisfy:\notka{kros-komut}
\begin{equation} \label{kros-komut} 
[A_{\kropka{c}}, \hat{f}]=\iota (\kot(X^r_{\kropka{c}})f)\,\widehat{}
\end{equation}
\end{enumerate}
\end{prop}\normalsize
\noindent{\em Proof:} 
Since $\widehat{u_{\kropka{c}}}$ is a one-parameter group of multipliers in $C^*_r(G_B)$ for its strict continuity it is sufficient to check continuity at $t=0$ of the mapping
$\R\ni t\mapsto \widehat{u_{\kropka{c}}}(t) (\omega)\in C^*_r(G_B)$ for $\omega \in \sA(G_B)$. Let us choose $\om_0=\lo\mt\ro$ as in (\ref{def-lo}, \ref{def-ro})\notka{def-lo,def-ro}, 
then $\om= f\om_0$  for $f\in\sD(G_B)$ and  we can write $\widehat{u_{\kropka{c}}}(t)(f\om_0)=:(\widehat{u_{\kropka{c}}}(t) f )\om_0$, where the function $(\widehat{u_{\kropka{c}}}(t) f)$ 
is given by (compare (\ref{app-Bc0-action})):
\notka{Bc0-action-1}
\begin{align}\label{Bc0-action-1} (\widehat{u_{\kropka{c}}}(t)f)( b c) & =f(u_{\kropka{c}}(-t) bc) \, j_C(\exp(t \kropka{c}))^{-1/2} 
\end{align}
Since the mapping $\R\times G_B\ni (t, g)\mapsto u_{\kropka{c}}(t) g\in G_B$ is continuous and $f\in \sD(G_B)$, for $\delta>0$ and $|t|<\delta$ supports of all functions
$(\widehat{u_{\kropka{c}}}(t)f)$ are contained in a fixed compact set, so by the lemma \ref{lemma-ind-lim}, it is sufficient to prove 
that $(\widehat{u_{\kropka{c}}}(t)f)$ converges, as $t\rightarrow 0$,  uniformly to $f$. But this is clear, since everything happens in a fixed compact set and all functions and mappings 
appearing in (\ref{Bc0-action-1}) are smooth.

Recall that the  domain of a generator $A$ of a (strongly continuous) one-parameter group of unitaries $U_t$ on a Hilbert space is defined as 
set of those vectors $\psi$ for which the limit $\displaystyle\lim_{t\rightarrow 0}(-\iota)\frac{U_t\psi-\psi}{t}$  exists and $A\psi$ is  the value of this limit.
Moreover if a dense linear subspace is contained in the domain of $A$ and is invariant for all $U_t$'s then it is a core for $A$.

Let us choose $\nu_0$ -- real non vanishing \halden  on $B$ and let $\Psi_0=\rho_0\mt\nu_0$. This is (real, non vanishing) \halden on $G_B$ 
and any $\psi\in \omh(G_B)$ can be written as $\psi=f \Psi_0$ for $f\in \sD(G_B)$.  The action of $\widehat{u_{\kropka{c}}}(t)$ can be written as 
$\widehat{u_{\kropka{c}}}(t)(f\Psi_0 )=:(\widehat{u_{\kropka{c}}}(t)f)\Psi_0$ and $(\widehat{u_{\kropka{c}}}(t)f)$ is given by the formula (\ref{app-Bc0-action}).
Using this formula and formulae (\ref{def-u-Xr-ckropka}, \ref{app-def-modular-functions}) one verifies that: \notka{A-formula-0}
\begin{align}\label{A-formula-0}
\lim_{t\rightarrow 0}\frac{1}{t}\left((\widehat{u_{\kropka{c}}}(t)f)( b c)- f(b c)\right) &= -\left((X^r_{\kropka{c}} f)( bc) + \frac{1}{2} Tr(ad(\kropka{c})|_{\gotc}) f(bc)\right),
\end{align}
where $ad(\kropka{c})(\kropka{e}):=[\kropka{c},\kropka{e}]\,,\,\kropka{c},\kropka{e}\in\gotc$. Again, since in the formula above, we stay, for a given $f$,  in a fixed compact 
subset and everything is smooth,  the limit  is, in fact, uniform and therefore also in $L^2(G_B)$.
Thus $\omh(G_B)\subset Dom(A_{\kropka{c}})$ and \notka{A-formula-1}
\begin{equation}\label{A-formula-1}
A_{\kropka{c}}(f\,\Psi_0)=:A_{\kropka{c}}(f) \,\Psi_0 \,,\,\,A_{\kropka{c}}(f):=\iota\left(X^r_{\kropka{c}} f + \frac{1}{2} Tr(ad(\kropka{c})|_{\gotc}) f\right)
\end{equation}
%%%%%%%%%%
Clearly $\omh(G_B)$ is dense in $L^2(G_B)$  and $\widehat{u_{\kropka{c}}}(t)$ invariant, so it is a core for $A_{\kropka{c}}$.\\
%%%
Since $\widehat{u_{\kropka{c}}}(t)(\psi)$ is a {\em push-forward} by the flow of $X^r_{\kropka{c}}$, 
the limit $\displaystyle \lim_{t\rightarrow 0}\frac{1}{t}\left((\widehat{u_{\kropka{c}}}(t)\psi)(g) -\psi(g)\right)$ 
is equal to $- \sL_{X^r_{\kropka{c}}} \psi$ and, consequently, $A_{\kropka{c}}\psi= \iota \sL_{X^r_{\kropka{c}}}\psi$.

Let us prove (\ref{komut-A}).  The mapping $C\times G_B\ni(c_0,g)\mapsto (B c_0) g\in G_B$ is a left action of $C$. For $\kropka{c}\in \gotc$, the fundamental vector 
field for this action is defined as $\displaystyle X_{\kropka{c}}(g):= \left.\frac{d}{d t}\right|_{t=0} u_{\kropka{c}}(-t)(g)$ and the map $\kropka{c}\mapsto X_{\kropka{c}}$
is a Lie algebra homomorphism (see e.g \cite{Lie-Mar})  Comparing with (\ref{def-u-Xr-ckropka}) we see that 
$X^r_{\kropka{c}}=(-1) X_{\kropka{c}}$. Since  $A_{\kropka{c}}=\iota \sL_{X^r_{\kropka{c}}}$ the formula (\ref{komut-A}) follows.

The formula (\ref{kros-komut}) is a direct consequence of (\ref{bis-fun-komut2}) and (\ref{def-kotwica}).\\
\dow
%}
%%%%%%%%%%%%%%%%%%%%%%%%%%%%%%%%%%%%%%%%%%%%%%%%%%%%%%%%%%%%%%%%%%%%%%%%%%%%%%%%%%%%%%%%%%%%%%%%%%%%%%%%%%%%%%%%%%%%%%
\subsection{Commutation relations for  $\kappa$-Poincar\'{e}}

Let us consider the one-parameter group $\,c(t):=\exp (t \Mlambda_{0(n+1)})$; in $(s,y)$ coordinates $c(t)=(s(t), y(t)):=(e^{-t},0)$.  Let $S(t)=B c(t)$ be 
the corresponding group of bisections i.e.
\begin{equation*}
S(t):=\{ (\Lambda, u, w,\alpha;e^{-t},0)\,,\, (\Lambda, u, w,\alpha)\in B\,,\,t\in \R\}.
\end{equation*}
%%%%%%%%%%%%%%%%%%%%%

\noindent
Define $\displaystyle (\tilde{\Lambda}(t),\tilde{u}(t),\tilde{w}(t),\tilde{\alpha}(t)):=b_R(\Lambda, u, w,\alpha;e^{t},0)$ and let $\hat{S}$ be  the generator of $\widehat{S(t)}$.
By (\ref{Bc0-action0},\ref{def-kotwica}) and (\ref{kros-komut}) we get:
$$[\hat{S},\hat{Q}](\Lambda,u, w ,\alpha)=\iota \left(\left.\frac{d}{d t}\right|_{t=0} \tilde{Q}(t)\right)^{\widehat{}}\,,\, Q=\Lambda,u, w,\alpha$$
By  (\ref{solution-iwasawa1}) we have 
%%%%%%%%%%%%%%%%%%%%%
\begin{align*}
\tilde{\Lambda}(t) &= \Lambda -\frac{\sinh t}{\cosh t  +\alpha \sinh t} u w^t & 
\tilde{\alpha}(t) &= \frac{\alpha\cosh t+\sinh t}{\cosh t+\alpha\sinh t} \\
\tilde{u}(t) &= \frac{u}{\cosh t+\alpha\sinh t} & \tilde{w}(t) &= \frac{w}{\cosh t+\alpha\sinh t},
\end{align*}
and by differentiation we obtain commutation relations (as operators on $\sA(G_B)$ or $\omh(G_B)$): \notka{relkom-S}
\begin{align}\label{relkom-S}
[\hat{S},\hat{\Lambda}] &= -\iota \widehat{u w^t} & [\hat{S},\hat{w}] &=  -\iota \hat{\alpha}  \hat{w} &   
   \,[\hat{S},\hat{u}] &= -\iota  \hat{\alpha}  \hat{u} &   [\hat{S},\hat{\alpha}] &=   -\iota  (\hat{\alpha}^2  -1)
\end{align}
or with indices put explicitly:\notka{relkom-S-1}
\begin{align}\label{relkom-S-1}
[\hat{S},\widehat{\Lambda_{kl}}] &= -\iota \widehat{u_k}\widehat{w_l} & [\hat{S},\widehat{w_k}] &=  -\iota \hat{\alpha}  \widehat{w_k} &   
   \,[\hat{S},\widehat{u_k}] &= -\iota  \hat{\alpha}  \widehat{u_k} &   [\hat{S},\hat{\alpha}] &=   -\iota  (\hat{\alpha}^2  -1)
\end{align}

Now, for  $y_0\in \R^n$, let $c(t):=\exp (t \Mlambda(y_0))$, where, as in (\ref{def-wsp-sy}), $\Mlambda(y_0):=\sum_{i=1}^n(y_0)_i (\Mlambda_{i(n+1)}-\Mlambda_{i 0})$);
or in $(s,y)$ coordinates: $c(t)=(s(t),y(t))= (1,t y_0)$.  Let  $ Y_0(t)$  be the  corresponding one-parameter group of bisections $B c(t)$; as before  we put 
$\displaystyle (\tilde{\Lambda}(t),\tilde{u}(t),\tilde{w}(t),\tilde{\alpha}(t)):=b_R(\Lambda, u, w, \alpha;1, -ty_0)$. 
By (\ref{solution-iwasawa1}) we get: 
%%%%%%%%%%%
\begin{align*}
\tilde{\Lambda}(t) &= \Lambda + \frac{t^2|y_0|^2}{2M(t)} u w^t - \frac{t(1+t w^t y_0)}{M(t)} u y_0^t -  \frac{t}{M(t)}\Lambda y_0 w^t- \frac{t^2(1-\alpha)}{M(t)} \Lambda y_0 y_0^t 
\end{align*}
\begin{align*}
\tilde{\alpha}(t) &= 1-\frac{1-\alpha}{M(t)}  &  \tilde{u}(t) &= \frac{u+ t( w^ty_0 u+ (1-\alpha)\Lambda y_0)}{M(t)} & \tilde{w} (t) &= \frac{w+ t (1-\alpha)y_0}{M(t)},
\end{align*}
$${\rm where}\,\,\,\,\,M(t):=\frac{|y_0|^2}{2}(1-\alpha)t^2 + t w^t y_0+1$$
Denoting by   $\widehat{Y_0}$  the generator of $\widehat{Y_0(t)}$ we obtain commutation relations (again as  operators on $\sA(G_B)$ or $\omh(G_B)$):
%%%%%%%%%%%%%%
\notka{relkom-Y}
\begin{align}\label{relkom-Y}\nonumber
[\widehat{Y_0},\hat{\Lambda}] &=  -\iota \left(u y_0^t+\Lambda y_0 w^t\right)^{\widehat{}}     &  [\widehat{Y_0},\hat{\alpha}] &= \iota\left ((1- \alpha)w^t y_0\right)^{\widehat{}}   \\
[\widehat{Y_0},\hat{w}] &= \iota\left((1-\alpha)y_0 -w^t y_0 w\right)^{\widehat{}}  & [\widehat{Y_0},\hat{u}] &= \iota \left((1- \alpha)\Lambda y_0\right)^{\widehat{}} 
\end{align}
%%%%%%%%%%%%%%%%%%%%%%%%%%%%%%%%%%%%%%%%%%%%%%%%%%%%%%%%%%%%%%%%%%%%%%%%%%%%%%%%%%%%%%%%%%%%%%%%%
or for $y_0:=e_m\in\R^n$ with corresponding $\widehat{Y_m}$  and with indices of matrix elements:\notka{relkom-Y-1}
\begin{align}\label{relkom-Y-1}\nonumber
[\widehat{Y_m},\widehat{\Lambda_{kl}}] &=  -\iota \left(\widehat{u_k} \delta_{ml}+\widehat{\Lambda_{km}}\widehat{w_l}\right)    &  
[\widehat{Y_m},\hat{\alpha}] &= \iota (1- \hat{\alpha})\widehat{w_m}  \\
[\widehat{Y_m},\widehat{w_k}] &= \iota\left((1-\hat{\alpha})\delta_{mk} -\widehat{w_m}\widehat{w_k}\right) & 
[\widehat{Y_m},\widehat{u_k}] &= \iota (1- \hat{\alpha})\widehat{\Lambda_{km}} 
\end{align}
%%%%%%%%%%%%%%%%%%%%%%%%%%%%%%%%%%%%%%%%%%%%%%%%%%%%%%%%%%%%%%%%%%%%%%%%%%%%%%%%%%%%%%%%%%%%%%%%

For completeness let us write relations between $\hat{S}$ and $\widehat{Y_0}$. Since $\hat{S}$ ($\widehat{Y_0}$) is the generator of the group $\widehat{u_{\kropka{c}}}$ for 
$\kropka{c}=\Mlambda_{0(n+1)}$ ($\kropka{c}=\Mlambda(y_0)$)  by the formula (\ref{komut-A}) and (\ref{lambda-komut}) we obtain:\notka{relkom-S-Y}
\begin{align}\label{relkom-S-Y}
[\hat{S}, \widehat{Y_0}] &= - i \widehat{Y_0}\,,& [\widehat{Y_1}, \widehat{Y_0}] &=0,
\end{align}
where $\widehat{Y_1}$ is defined in the same way as $\widehat{Y_0}$ for a vector $y_1\in\R^n$.

Our generators $\hat{S}$ and $\hat{Y}_k$ are related, via the mapping $\gotc\ni \kropka{c}\mapsto A_{\kropka{c}}$ used in the previous subsection, to the following basis in  $\gotc$:
\notka{def-generator-basis}
\begin{equation} \label{def-generator-basis}
(\hat{S}, \hat{Y}_k)\leftrightsquigarrow (\kropka{c}_0,\kropka{c}_k) :=(\Mlambda_{0(n+1)},\Mlambda_{k(n+1)}-\Mlambda_{k 0})\,,\quad k=1, \dots, n
\end{equation}
\subsection{Formulae for  comultiplication}

Now we consider again a general DLG $(G;B,C)$ together with the relation $\delta_0:G_B\rel G_B\times G_B$ defined in 
(\ref{def-delta0}) i.e. 
$\displaystyle \delta_0=\{(b_1^{-1}c_L(b_1 b c), b_1 b c; b c): b,b_1\in B, c\in C\}$. It defines the mapping $\hat{\delta}_0$ given by the formula
(\ref{app-hat-delta0-form}), which extends to the coassociative $\Delta_0\in Mor(C^*_r(G_B), C^*_r(G_B)\mt C^*_r(G_B))$ \cite{PS-DLG}.
For a smooth, bounded function $f$ on $B$, let $\Delta_B(f)$ be the value of the comultiplication of the group $B$ on $f$ i.e. $\Delta_B(f) (b_1,b_2):=f(b_1 b_2)$. 
Let $\hat{f}$ be the multiplier of $C^*_r(G_B)$ defined by $f$;  the formula (\ref{app-hat-delta0-form}) implies \notka{delta0-function}
\begin{align}\label{delta0-function}
\Delta_0(\hat{f})&=\widehat{\Delta_B(f)}
\end{align}

One easily computes the image of a bisection $B c_0$ by $\delta_0$:  
\begin{equation*}
\delta_0(B c_0)=\{(b_1c_L(b c_0), b c_0) : b, b_1\in B\}
\end{equation*}
and its action on $G_B\times G_B$: \notka{delta0-bis}
\begin{equation}\label{delta0-bis}
\begin{split}
\delta_0(B c_0)(b_1 c_1, b_2 c_2) & =(c_L^{-1}(b_1 c_L(b_2 c_0^{-1})) b_1 c_1, c_L^{-1}(b_2 c_0^{-1}) b_2 c_2)=\\
& = (b_R(b_1 c_L(b_2 c_0^{-1}) )c_L^{-1}(b_2 c_0^{-1}) c_1, b_R(b_2c_0^{-1}) c_0 c_2)=\\
& = (c_R(c_R(c_0b_2^{-1}) b_1^{-1}) b_1 c_1, c_R(c_0b_2^{-1})  b_2 c_2)
\end{split}
\end{equation}
%%%%%%%%%%%%%%%%%%%%%%
For  $\kropka{c}\in\gotc$,  let $u_{\kropka{c}}$ be  the one-parameter group of bisections defined in (\ref{def-u-Xr-ckropka}) 
and define $\tilde{u}_{\kropka{c}}(t):=\delta_0(u_{\kropka{c}}(t))$; 
this is a one-parameter group of bisections of $G_B\times G_B$. Let  $\tilde{X}^r_{\kropka{c}}$ be the corresponding  right invariant vector field on $G_B\times G_B$ 
(compare (\ref{def-u-Xr-ckropka})). By (\ref{delta0-bis}):
$$\tilde{X}^r_{\kropka{c}}(b_1 c_1, b_2 c_2)=\left.\frac{d}{d t}\right|_{t=0} (\tilde{c}_1(t) b_1 c_1,\tilde{c}_2(t) b_2 c_2),\, {\rm where} $$
$$\tilde{c}_1(t) := c_R( b_1 c_R( b_2 \exp( t\kropka{c}) b_2^{-1}) b_1^{-1}) \,,\,\,\tilde{c}_2(t) =c_R(b_2\exp(t\kropka{c})b_2^{-1}),$$
and  \notka{delta0-Xr}
\begin{equation}\label{delta0-Xr}
\tilde{X}^r_{\kropka{c}}(b_1 c_1, b_2 c_2)=(\adc(b_1)\adc(b_2)(\kropka{c}) b_1 c_1,\,\adc(b_2)(\kropka{c}) b_2 c_2)
\end{equation}
%%%%%%%%%%%%%%%%%%%%%%
Let $(X_\alpha)$  be a basis in $\gotc$ , $u_\alpha,\, A_\alpha,\, X^r_\alpha$  be corresponding one-parameter groups, 
their  generators and right-invariant vector fields respectively,  and let $\tilde{u}_\alpha:=\delta_0(u_\alpha)$ with corresponding  $\tilde{A}_\alpha$ and $\tilde{X}^r_\alpha$. 
The formula (\ref{delta0-Xr}) reads (after identification $T_{(x,y)}(X\times Y)=T_xX\oplus T_yY$) : \notka{delta0-Xr-1}
\begin{equation}\label{delta0-Xr-1}
\tilde{X}^r_\alpha(b_1 c_1, b_2 c_2)=\sum_\beta \adc_{\beta \alpha}(b_2)  X^r_\beta( b_1 c_1) + X^r_\alpha( b_2 c_2),
\end{equation}
where functions $ \adc_{\beta \alpha}: B\rightarrow \R$ are matrix elements of $\adc|_B$ (\ref{adb-adc}),  i.e.  $\adc(b)(X_\alpha)=:\sum_\beta \adc_{\beta \alpha}(b) X_\beta$. 
As in Prop.\ref{prop-generators} the action of $\tilde{A}_\alpha$ on 
$\omh(G_B\times G_B)$ is given by ($\iota$ times) the Lie derivative with respect to $\tilde{X}^r_\alpha$. Since $\tilde{A}_\alpha=\Delta_0(A_{\alpha})$ we obtain the following equality on
$\omh(G_B\times G_B)$  \notka{Delta0-gen}
\begin{equation}\label{Delta0-gen}
\Delta_0(A_{\alpha})=I\mt A_\alpha+ \sum_\beta A_\beta \mt \adc_{\beta \alpha}
\end{equation}
%%%%%%%%%%%%%%%%%%%
\begin{re}
Let us comment on the meaning of the equality above. The operator $\tilde{A}_\alpha=\Delta_0(A_{\alpha})$ is essentially self-adjoint on $\omh(G_B\times G_B)$. The
operator on the right hand side has immediate meaning on $\omh(G_B)\mt\omh(G_B)$ and is symmetric there.  But since elements of  $\,\omh(G_B\times G_B)$ can be approximated 
by elements of $\omh(G_B)\mt\omh(G_B)$  in topology of $\omh(G_B\times G_B)$ (i.e. uniformly with all derivatives on compact sets)  and operators $A_\alpha$ are differential operators, 
the space $\omh(G_B\times G_B)$ is in the  domain of the closure of the right hand side (treated as operator on $\,\omh(G_B)\mt\omh(G_B)$). 
Therefore this closure is the self-adjoint operator $\Delta_0(A_{\alpha})$.
\end{re}

\subsection{Comultiplication for $\kappa$-Poincar\'e}

For a bisection $B c_0$ let us consider the map: $(b_1 c_1,b_2 c_2)\mapsto \delta(B c_0)(b_1 c_1,b_2 c_2)$, where the relation $\delta$ is given by  (\ref{def-delta}). 
By the use of  that formula one finds  the domain of this map -- $\{(b_1 c_1 ,b_2 c_2): b_2, b_R(b_2 c_0^{-1})\in B'\}$ and  
\notka{delta-bis}
\begin{equation}\label{delta-bis}
\begin{split}
\delta(B c_0) (b_1 c_1,b_2 c_2)=&
\left(b_R[b_1 \tilde{c}_L^{-1}(b_2) \tilde{c}_L(b_2 c_0^{-1})] 
\tilde{c}_L^{-1}(b_2 c_0^{-1})\tilde{c}_L(b_2) c_1,\, b_R(b_2 c_0^{-1})c_0 c_2\right)=\\
=&\left( c_L^{-1}(b_1 \tilde{c}_L^{-1}(b_2) \tilde{c}_L(b_2 c_0^{-1}))b_1 c_1,\,c_L^{-1}(b_2 c_0^{-1})b_2 c_2 \right).
\end{split}
\end{equation}
Note that $b_R(b_2 c_0^{-1})c_0 c_2=(B c_0)(b_2 c_2)$ i.e. the action
on the right ``leg'' is the action of the bisection $B c_0$.

We will use this expression to get the  comultiplication for generators. Let  $c_0:=c(t):=\exp(t \kropka{c})$ for  $\kropka{c}\in\gotc$.
Since $B'$ is open in $B$, for fixed $b_2\in B'$ the right hand side of (\ref{delta-bis}) is well defined for $t$ sufficiently close to $0$; therefore we can define 
a vector field on $G\times B'C$ which we denote by 
$\delta(X^r_{\kropka{c}})$. It is given by:\notka{delta-Xr}
\begin{align*}\label{delta-Xr}
\delta(X^r_{\kropka{c}})(b_1 c_1, b_2 c_2) &:=\left.\frac{d}{d t}\right|_{t=0} (c_1(t) b_1 c_1,c_2(t)b_2 c_2), & \\\nonumber
c_1(t)& :=c_L^{-1}(b_1\tilde{c}_L^{-1}(b_2)\tilde{c}_L(b_2\exp(-t \kropka{c})))\,,&  c_2(t) &:=c_L^{-1}(b_2 \exp(-t \kropka{c}))\end{align*}
Let $\trzutC$ be the projection onto $\gotc$ corresponding to the decomposition $\gotg=\gotc\oplus\gota$ and \notka{def-tadc}
\begin{equation}\label{def-tadc}
\tadc(g):=\trzutC Ad(g)|_{\gotc}.
\end{equation}
Computing  derivatives we obtain
\notka{delta-Xr-1}
\begin{equation}\label{delta-Xr-1}
\delta(X^r_{\kropka{c}})(b_1 c_1, b_2 c_2)=\left((\adc(b_1)\tadc(a_R(b_2) \kropka{c}) b_1 c_1\,, \, (\adc(b_2)\kropka{c}) b_2 c_2\right).
\end{equation}
For a basis $(X_\alpha)$ in $\gotc$ we obtain (compare \ref{delta0-Xr-1}): \notka{delta-Xr-2}
\begin{equation}\label{delta-Xr-2}
\delta(X^r_\alpha)(b_1 c_1, b_2 c_2)=\sum_\beta \tadc_{\beta \alpha}(a_R(b_2)) X^r_\beta( b_1 c_1) + X^r_\alpha( b_2 c_2),
\end{equation}
%%%%%%%%%%%%%%%%%%%%%%%%%%%

We may try to  go one step further and   write this equation in a form similar to (\ref{Delta0-gen}) as: \notka{Delta-gen}
\begin{equation}\label{delta-gen}
\Delta(A_\alpha) = I\mt A_\alpha + \sum_\beta A_\beta\mt (\tadc_{\beta \alpha} \cdot a_R)
\end{equation}
%%%%%%%%%%%%
But, contrary to (\ref{Delta0-gen}), where the right hand side has a well defined meaning, here we have rather formal expression. 
Certainly we have equality as operators on $\omh(G_B)\mt\omh(B'C)$ but {\em this space  is not a core for the left hand side} so it is not 
true that the closure of the right hand side is equal to the self-adjoint operator on the left hand side.
The precise formula for  comultiplication is given in Prop.\ref{prop-Delta}.

Let us now, compute the  comultiplication for our generators  $(\widehat{S},\widehat{Y}_k)$ (formally, in the sense of (\ref{delta-gen})). Recall that they are related (\ref{def-generator-basis}) to the following basis in $\gotc$:
\begin{equation*}%\label{def-generator-basis}
(\kropka{c}_\beta):=(\kropka{c}_0,\kropka{c_k}):=(\Mlambda_{0(n+1)},\Mlambda_{k(n+1)}-\Mlambda_{k 0})\,,\quad k=1, \dots, n.
\end{equation*}
%%%%%%%%%%%%%%%%
We need to find matrix elements of representation $\tadc(a),\,a\in A$ in the basis $(\kropka{c}_\beta)$. Let us denote this matrix by $\kkad(a)$:
$$\tadc(a)\kropka{c}_\beta=\sum_\alpha \kkad_{\alpha\beta}(a) \kropka{c}_\alpha$$
%%%%%%%%%%%%%%%%
We will use lemma \ref{lema-reps}.
It is easy to see that the orthogonal complement of $\gota$, with respect to the form $k$ defined in (\ref{app-def-k}), is 
\mbox{$\displaystyle \gota^\perp=span \{\Mlambda_{\beta(n+1)}\,:\,\beta=0,\dots, n\}$} and bases $(e_\beta)$ and $(\rho_\beta)$ defined as:
\notka{def-bases-e-ro}
\begin{equation}\label{def-bases-e-ro}
e_\beta:=\Mlambda_{\beta(n+1)}\quad,\quad\quad\rho_\beta:=k(e_\beta)\quad,\quad\quad \beta=0,\dots, n
\end{equation}
are  orthonormal basis in $\gota^\perp$ and $\gota^0$, respectively.
The projection  $\trzutC:\gota^\perp \rightarrow \gotc$ acts as  $\trzutC(e_\beta)=\Mlambda_{\beta (n+1)} -\Mlambda_{\beta 0}=\kropka{c}_\beta$.
By the lemma \ref{lema-reps}, matrix of $\tadc(a)$ (i.e. the matrix $\kkad(a)$) is equal to the matrix of $\kad(a)|_{\gota^0}$ in basis $(\rho_\alpha)$.
%%%%%%%%%%%%%%%
For $a=(z,U,d)\in A$, by direct computations (e.g. using formulae in the Appendix) one gets:
\notka{kad-zUd}
\begin{align}\label{kad-zUd}
\kkad(z,U,d)=\left(\begin{array}{ccc}
d \frac{1+|z|^2}{1-|z|^2} & \frac{2}{1-|z|^2}z^tU D_1 \\ 
d \frac{2}{1-|z|^2} z  & (I+\frac{2}{1-|z|^2}z z^t)U D_1 \end{array}\right),\quad 
{\rm where}\quad D_1=\left(\begin{array}{cc}d I_{n-1} & 0\\0 & 1\end{array}\right).
\end{align}
Finally, using (\ref{def-aR}) one gets:\notka{kad-zUd-ar}
\begin{equation}\label{kad-zUd-ar}
\kkad(a_R(\Lambda,u,w,\alpha))=\left(\begin{array}{ccc}
\frac{1}{\alpha} & \frac{w^t}{\alpha} \\ 
- \frac{u}{|\alpha|}  & sgn(\alpha)( \Lambda-\frac{u w^t}{\alpha}) \end{array}\right)
\end{equation}
%%%%%%%%%%%%%%%%
And formulae for comultiplication on generators are:
%%%%%%%%%%%%%%
\notka{Delta-na-generatorach}
\begin{align}\label{Delta-na-generatorach}
\begin{split}
\Delta(\widehat{S}) &= I\mt \widehat{S}+ \widehat{S}\mt \frac1\alpha+\sum_k\widehat{Y}_k\mt \frac{- u_k}{|\alpha|}\\
\Delta(\widehat{Y}_i) &= I\mt\widehat{Y}_i+\widehat{S}\mt \frac{w_i}{\alpha} + \sum_k\widehat{Y}_k\mt sgn(\alpha)(\Lambda_{ki}-\frac{u_k w_i}{\alpha})
\end{split}
\end{align}
%%%%%%%%%%%%%%%%%%%%%%%%%
It remains to compute $\Delta$ for generators $(\Lambda,u,w,\alpha)$. Let $f$ be a  smooth and bounded function on $B$. By  (\ref{def-delta}) $\Delta(f)$, as a function on $B\times B'$, 
is given by the formula: \notka{Delta-functions-B}
\begin{equation}\label{Delta-functions-B}
(\Delta f) (b_1, b_2)=f(b_R(b_1 a_R(b_2)))=f(b_R(b_1(\tilde{c}_L(b_2))^{-1}) b_2)\,,\,\,(b_1,  b_2)\in B\times B'
\end{equation}
%%%%%%%%%%%%%%
Using formulae (\ref{def-tcL}) and (\ref{solution-iwasawa1}), after some computations, one gets:\notka{Delta-matrix-elem-B}
\begin{align}\label{Delta-matrix-elem-B}
\Delta(u_k) &= u_k\mt\sign(\alpha)+ P^{-1} \sum_l \alpha\Lambda_{kl}\mt u_l\,,\quad\Delta(w_k) =I\mt w_k+ P^{-1} \sum_l w_l\mt |\alpha| \Lambda_{lk}  & \nonumber\\
\Delta(\alpha) &=P^{-1}(\alpha\mt\alpha)\,,  \quad \quad\Delta(\Lambda_{kl}) =\sum_j \Lambda_{kj}\mt\Lambda_{jl} + P^{-1}\sum_{mj}\Lambda_{km} w_j\mt \sign(\alpha)u_m\Lambda_{jl}\,, &
\end{align}
where $\displaystyle P:=1-\sum_k w_k\mt \sign(\alpha)u_k$. Notice that $P$ is invertble on $B\times B'$ and right hand sides are well defined (as smooth functions on $B\times B'$).
%%%%%%%%%%%%%%%%%%%%%%%%%%%%%%%%%%%%%%%%%%%%%%%%%%%%%%%%%%%%%%%%%%%%%%%%%%%%%%%%%%%%%%%%
%%%%%%%%%%%%%%%%%%%%%%%%%%%%%%%%%%%%%%%%%%%%%%%%%%%%%%%%%%%%%%%%%%%%%%%%%%%%%%%%%%%%%%%
\subsection{A closer look at the twist}
%%%%%%%%%%%%%%%%%%%%%%%%%%%%%%%%%%%%%%%%%%%%
\newcommand{\cexp}{\exp_C}
\newcommand{\clog}{\log_C}
%%%%%%%%%%%%%%%%%%%%%%%%%%%%%%%%%%%%%%%%%%%%%%%%%
\begin{center} {\em In this subsection groups  $G,A,B,C$ fulfill  assumptions (\ref{basic-assumpt}) and  $C$ is an exponential Lie group i.e. 
the exponential mapping $\cexp: \gotc \rightarrow C$ is a diffeomorphism.}
\end{center}
Let $\clog:=\exp_C^{-1}$ and for  $b\in B'$ let \notka{def-ct}
\begin{equation}\label{ct-def}
c_t(b):=\cexp(t \clog(\tilde{c}_L(b)^{-1}))=\cexp(- t \clog(\tilde{c}_L(b)))\,,\quad t\in\R
\end{equation}
be the one-parameter group defined by $-\clog(\tilde{c}_L(b))\in \gotc$ and 
%%%%%%%%%%%
\notka{def-Tt}
\begin{equation} \label{def-Tt}
T_t:=\{(b_1 c_t(b_2), b_2):b_1\in B, b_2\in B'\}\subset G_B\times\Gamma_{B'} .
\end{equation}
Note that $T_1$ is the twist (\ref{def-twist}).
%%%%%%%%%%%%%%%%%%%%%%%%%%%%%%%%%%%%%%%%%%%%%%%%%%%%%%%%%%
\begin{lem}
$T_t$ is a one-parameter group of bisections  $G_B\times\Gamma_{B'}$.
\end{lem}
{\em Proof:} It is easy to verify that  $T_t$ is a bisection and a submanifold for any $t\in\R$. For $s,t\in \R$:
$$(b_1 c_1, b_2 c_2)\in T_t T_s \iff \exists\, b_3,b_5\in B\,, b_4, b_6\in B': (b_1 c_1, b_2 c_2)=(b_3 c_t(b_4), b_4) (b_5 c_s(b_6), b_6),$$
(on the right hand side there is the multiplication in $G_B\times \Gamma_{B'}$)
so $c_2=e$ and $b_2=b_4=b_6\in B'$, and $b_1 c_1= m_B(b_3 c_t(b_6),b_5 c_s(b_6))$. Therefore $b_5=b_R(b_3 c_t(b_6))$ and $b_1 c_1=b_5 c_t(b_6)c_s(b_6)=b_5 c_{t+s}(b_6)$, i.e.
$T_t T_s =\{(b_1 c_{t+s}(b_2), b_2): b_1\in B, b_2\in B'\}=T_{t+s}$. \\
\dowl
%%%%%%%%%%%%%%%%%%%%%%%%%%%%%%%%%%%%%%%%%%%%%%%%%%%%%%%%%%%%%%%

\noindent
Let $X_T^r$ be the right invariant vector field on $G_B\times\Gamma_{B'}$ defined by $T_t$ (compare (\ref{def-u-Xr-ckropka})) i.e. 
$$X_T^r(b_1 c_1, b_2 c_2):=\left.\frac{d}{d t}\right|_{t=0} T_t(b_1 c_1, b_2 c_2)$$
Using the definition (\ref{def-Tt}) we get
for  $b_2\in B'$:
\begin{align*}\begin{split}
T_t(b_1 c_1, b_2 c_2) & = (b_R(b_1c_t(b_2)^{-1})c_t(b_2) c_1, b_2 c_2)=(c_L(b_1 c_t(b_2)^{-1})^{-1} b_1 c_1\,,\,b_2 c_2)=\\
&=( c_R(c_t(b_2)b_1^{-1}) b_1 c_1\,,\,b_2 c_2), \,\,
\end{split}
\end{align*}
and \notka{twist-pole-0}
\begin{align}\label{twist-pole-0}
X_T^r(b_1 c_1, b_2 c_2) & =\left((-\adc(b_1) \log_C(\tilde{c}_L(b_2))) b_1 c_1, 0_{b_2 c_2}\right)
\end{align}
For a basis $(X_\alpha)$ of $\gotc$ let
% $D_\alpha: \G_{B'}\rightarrow  \R$   RACZEJ TAK
$D_\alpha: B'\rightarrow  \R$  be coordinates of $\clog(\tilde{c}_L(b))$:\notka{def-D}
\begin{equation}\label{def-D}
\clog(\tilde{c}_L(b))=:\sum_\alpha D_\alpha(b) X_\alpha.
\end{equation}
Using the definition (\ref{def-u-Xr-ckropka}) of $X^r_\alpha$ we can write (\ref{twist-pole-0}) as \notka{twist-pole}
\begin{equation}\label{twist-pole}
X_T^r(b_1 c_1, b_2 c_2)=-\sum_\alpha D_\alpha (b_2) X^r_\alpha(b_1 c_1)
\end{equation}
\begin{prop}
Let $\widehat{T}_t$ be the one-parameter group of unitaries in 
$L^2(G_B\times \Gamma_{B'})=L^2(G_B)\mt L^2(\Gamma_{B'})=L^2(G_B)\mt L^2(G_B$) defined by $T_t$. The action of $\widehat{T}_t$ on $\sA(G_B\times \Gamma_{B'})$ is given by
\begin{align}\widehat{T_t}(F(\om_0\mt\om_0))&=:(\widehat{T_t}F)(\om_0\mt\om_0)\,,& (\widehat{T_t}F)(g_1,g_2)& := F(T_{-t}(g_1,g_2))\, j_C(c{-t}(b_L(g_2)))^{-1/2},
\end{align}
where $\om_0=\lambda_0\mt\rho_0$ for  $\lambda_0$ and $\rho_0$  defined in  (\ref{def-lo},\ref{def-ro}) and $F\in \sD(G_B\times \Gamma_{B'})$.\\
%\tilde{c}_t:=\cexp(t\clog\tilde{c}_L(b_L(g_2)))$$
Let $\sT$ be the  generator of  $\widehat{T}_t$ .  $\sT$ ie essentially self-adjoint on $\omh(G_B\times \Gamma_{B'})$  and for $F(\psi_0\mt\psi_0) \in \omh(G_B\times \Gamma_{B'})$:
\notka{twist-generator}
\begin{equation}\label{twist-generator}
\begin{split}
\sT (F(\psi_0\mt\psi_0))& =:(\sT F)(\psi_0\mt\psi_0)\,,\\
(\sT F) (b_1 c_1, b_2 c_2) &=\iota\left( (X_T^r F)(b_1c_1, b_2 c_2)-\frac{1}{2}Tr (ad(\clog(\tilde{c}_L(b_2)))|_{\gotc})F(b_1 c_1, b_2 c_2)\right)
\end{split}
\end{equation}
where $\psi_0=\rho_0\mt\nu_0$ for some real, non vanishing \halden $\nu_0$ on $B$. Moreover on $\omh(G_B\times \Gamma_{B'})$ \notka{twist-generator-1}
\begin{equation}\label{twist-generator-1}
\sT=- \sum_\alpha A_\alpha\mt D_\alpha
\end{equation}
\end{prop}
{\em Proof:} 

As for any bisection   $\widehat{T_t}F$ is given by \cite{DG} :
$$(\widehat{T_t}F)(T_t(g_1,g_2))=F(g_1,g_2)\frac{(\rho_0\mt\rho_0) (v g_1\mt w
  g_2)}{(\rho_0\mt\rho_0)(T_t(v g_1\mt w g_2))}\,,
\,v,w\in\Lambda^{max}(T_eC)$$

Let $c(s) \subset C$ be a curve with $c(0)=e$. For $(g_1,g_2)\in G_B\times \Gamma_{B'}$ let $(\tilde{g}_1 , g_2):=T_t(g_1, g_2)$. By (\ref{def-Tt}) we get
$$T_t(c(s) g_1, g_2)=(c_1(s) \tilde{g}_1, g_2)\quad\,\quad\quad  T_t( g_1, c(s) g_2)=(c_2(s) \tilde{g}_1 , c(s) g_2),$$
for some curves $c_1(s), c_2(s)$.
Identifying tangent spaces to right fibers with $T_eC$ in corresponding points,
one sees that the map we have to consider has the form 
$\left(\begin{array}{cc} M_1 & M_2\\0 & I\end{array}\right)$, so its
action on densities is determined by $M_1$ i.e. (derivative of)
$c(s)\mapsto c_1(s)$. But this is exactly as the action of the  bisection $Bc_t(b_L(g_2))$ (compare (\ref{Bc0-action0}) and by (\ref{app-Bc0-action}) we obtain:
$$(\widehat{T_t}F)(b_1 c_1 ,b_2 c_2)= F(T_{-t}(b_1 c_1, b_2 c_2))\, j_C(c_{t}(b_2))^{-1/2}$$
Since $c_t(b_2)=\cexp(- t \clog(\tilde{c}_L(b_2)))$ we can proceed exactly as in the proof of Prop.\ref{prop-generators} 
and by putting $\kropka{c}=-\clog(\tilde{c}_L(b_2))$ in formula (\ref{A-formula-1}) we obtain (\ref{twist-generator}). As before, 
since $\omh(G_B\times \Gamma_{B'})$ is invariant for $\widehat{T}_t$ it is a core for $\sT$. 
Formula (\ref{twist-generator-1}) is a direct consequence of (\ref{twist-pole}) and (\ref{twist-generator}).\\
\dowl
%%%%%%%%%%%%%%%%%%%%%%%%%%%%%%%%%%%%%%%%%%%%%%%%%%%%%%%%%%%%%%%%%%%%%%%%%%%%%%%%%%%%%%%%%%%5
\subsection{Twist for $\kappa$-Poincar\'e}

Let $(\kropka{s},\kropka{y})\in \R\times \R^n$ denotes element of $\gotc$. The exponential mapping $\cexp$ and its inverse $\clog$ are given by:
\begin{equation*}
\cexp: \gotc\ni(\dot{s},\dot{y})\mapsto (\exp(\dot{s}),\dot{y}\frac{\exp(\dot{s})-1}{\dot{s}})\in C\quad,\quad
\clog: C\ni  (s,y)\mapsto (\log(s),y\frac{\log(s)}{s-1})\in \gotg
\end{equation*}
With respect to the basis (\ref{def-generator-basis}) the decomposition is:
$\displaystyle (\kropka{s}, \kropka{y})=-\kropka{s}\kropka{c}_0 + \sum_k\kropka{y}_k\kropka{c}_k$. %(\Mlambda_{k(n+1)}-\Mlambda_{k 0})$
%%%
Recall  from (\ref{def-tcL}) the mapping $\tilde{c}_L: B'\ni (\Lambda,u,w,\alpha) \mapsto (|\alpha|,-sgn(\alpha) u )\in C$.
Thus \notka{D-for-kappa}
\begin{align}\label{D-for-kappa}
\clog(\tilde{c}_L(\Lambda,u,w,\alpha)) &= (\log(|\alpha|), -sgn(\alpha) \frac{\log(|\alpha|)}{|\alpha|-1} u) 
=  - \log(|\alpha|)\kropka{c}_0 - sgn(\alpha) \frac{\log(|\alpha|)}{|\alpha|-1}\sum_k u_k \kropka{c}_k \nonumber \\
 {\rm i.e.}\quad\quad & D_0(\Lambda,u,w,\alpha) =-\log(|\alpha|)\,,\,\,\,D_k(\Lambda,u,w,\alpha) = - sgn(\alpha) \frac{\log(|\alpha|)}{|\alpha|-1} u_k
\end{align}
The bisection $T_t$ is equal:
$$T_t:=\{(b_1, |\alpha_2|^{-t}, \frac{\sign(\alpha_2)(|\alpha_2|^{-t}-1)}{1-|\alpha_2|} u_2;
\Lambda_2,u_2,w_2,\alpha_2,1,0):\alpha_2\neq 0\,,\,b_1\in B\}\subset G_B\times\Gamma_{B'}$$
And from (\ref{twist-generator-1}) and (\ref{D-for-kappa}) we get the formula for the twist:
$$\hat{T}=\exp(\iota \sT)\,,\quad \sT=\widehat{S}\mt\log|\alpha|+\sum_k\widehat{Y}_k\mt \frac{\sign(\alpha) \log|\alpha|}{|\alpha|-1}u_k$$
%%%%%%%%%%%%%%%%%%%%%%%%%%%%%%%%%%%%%%%%%%%%%%%%%%%%%%%%%%%%%%%%%%%%%%%%%%%%%%%%%%%%%%%%%%%%%%%%%%%%%%%%%%%%%%%%%%%%%%%%%%%5
%%%%%%%%%%%%%%%%%%%%%%%%%%%%%%%%%%%%%%%%%%%%%%%%%%%%%%%%%%%%%%%%%%%%%%%%%%%%%%%%%%%%%%%%%%%%%%%%%%%%%%%%%%%%%%%%%%%%%%%%%%
%%%%%%%%%%%%%%%%%%%%%%%%%%%%%%%%%%%%%%%%%%%%%%%%%%%%%%%%%%%%%%%%%%%%%%%%%%%%%%%%%%%%%%%%%%%%%%%%%%%%%%%%%%%%%%%%%%%%%%%%%
\section{Comparison with  relations coming from Poisson-Poincar\'{e} Group}

The classical Poincar\'{e} Group, our $(C^*(G_B),\Delta)$ is a quantization of, is $(TA)^0\subset T^*G$. After identifying $T^*G$ with the semidirect product $\gotg^*\rtimes G$
(\ref{app-semi-direct}), $(TA)^0$ is identified with $\gota^0\rtimes A$. Invariant, bilinear form $k$ (\ref{app-def-k}) on $\gotg$ induces the form on $\gota^0$ which makes it 
a vector Minkowski space and the action of $A$ is  by orthogonal transformations. The bundle $(TA)^0$ is dual to the Lie algebroid of groupoid $\Gamma_A$, 
so objects related to this groupoid have more direct relation to  functions on $\gota^0\rtimes A$ (see \cite{PS-poisson} for a detailed description).  
The groupoid $G_B$ was used to overcome some functional analytic problems, let us now transport  our expressions  to the  groupoid $\Gamma_A$ to relate them 
to  formulae in \cite{SZ-94, Kosinski-Maslanka}.

\subsection{Back to $\Gamma_A$ picture}
{\em Right invariant fields $X^r_{\kropka{c}}$.}   By (\ref{Bc0-action0}) and (\ref{def-u-Xr-ckropka}) 
$X^r_{\kropka{c}}(bc)=\left.\frac{d}{d t}\right|_{t=0} b_R(b c_t^{-1})c_tc=\left.\frac{d}{d t}\right|_{t=0} c_L^{-1}(b c_t^{-1})b c$, where 
$c_t:=\exp(t \kropka{c}).$  If $bc\in\Gamma_{B'}$ i.e. $b, b_R(bc)\in B'$ then it is easy to see that $b_R(bc_t^{-1})\in B'$ for $t$ sufficiently close to $0$ and, consequently,
$b_R(b c_t^{-1})c_tc \in \Gamma_{B'}$. So we can push 
$X^r_{\kropka{c}}(bc)$ to $\Gamma_A$ by the mapping (\ref{wzor-embed}) $bc\mapsto a_R(b) c$: 
$$b_R(bc_t^{-1})c_t c  \mapsto a_R(b_R(bc_t^{-1})) c_t c=a_R(a_R(b)c_t^{-1}) (c_t a_R(b)^{-1}) a_R(b)c= \tilde{c}^{-1}_L(a_R(b)c_t^{-1}) a_R(b)c$$
since  $a_R(b_R(bc_t^{-1}))=a_R(bc_t^{-1})=a_R(a_R(b)c_t^{-1})$.
Let us denote the resulting (right-invariant) vector field on $\Gamma_A$ by $X^{A,r}_{\kropka{c}}$:\notka{def-Xr-kappa}
\begin{align}\label{def-Xr-kappa}
X^{A,r}_{\kropka{c}}(ac) &:=\left.\frac{d}{d t}\right|_{t=0} \tilde{c}^{-1}_L(a \exp(-t \kropka{c})) a c=(\widetilde{Ad^{\gotc}}(a)\kropka{c}) ac,
\end{align}
where $\tadc$ was defined in (\ref{def-tadc}). 

{\em Anchor.} Let us transfer formula (\ref{def-kotwica}): \notka{def-kotwica-kappa}
\begin{equation*} \label{def-kotwica-kappa}
\begin{split}
\left.\frac{d}{d t}\right|_{t=0} a_R(b_R(bc_t^{-1}))= \left.\frac{d}{d t}\right|_{t=0} a_R(a_R(b)c_t^{-1}) &=\left.\frac{d}{d t}\right|_{t=0} a_R(a_R(b)c_t^{-1}a^{-1}_R(b))a_R(b)\\
& =(-Ad^{\gota}(a_R(b))\kropka{c})a_R(b)
\end{split}
\end{equation*}
Thus, denoting the anchor of $\Gamma_A$ by $\kot^A$, we obtain: \notka{def-kotwica-kappa-1}
\begin{equation} \label{def-kotwica-kappa-1}
\kot^A(X^{A,r}_{\kropka{c}})(a):=(-Ad^{\gota}(a)\kropka{c})a.
\end{equation}

Comparing this formula and (\ref{def-Xr-kappa}) with (\ref{kotwica1}) and (\ref{def-u-Xr-ckropka}) we see that we can immediately write commutation relations analogous to relations
(\ref{komut-A}) and (\ref{kros-komut}): \notka{rel-komut-Gamma-A}
\begin{align} \label{rel-komut-Gamma-A}
[A^A_{\kropka{c}}, A^A_{\kropka{e}}] & =-\iota A^A_{[\kropka{c},\kropka{e}]},\,\kropka{c},\kropka{e}\in\gotc&\, 
[A^A_{\kropka{c}}, \hat{f}] &=\iota (\kot^A(X^{A,r}_{\kropka{c}})f)\,\widehat{}\,,\,f\in C^\infty(A).
\end{align}
Operators $A^A_{\kropka{c}}$ are, as before,  Lie derivatives along $X^{A,r}_{\kropka{c}}$ multiplied by $i$.  
These relations hold on $\omh(\Gamma_A)$ (and $\sA(\Gamma_A)$), however, their status is ``weaker'' then  one of relations (\ref{komut-A}) and (\ref{kros-komut}.
Firstly, operators $A^A_{\kropka{c}}$ {\em are not essentially self-adjoint} on $\omh(\Gamma_A)$; secondly, even though each smooth function on $A$ defines an element affiliated to 
 $C^*_r(\Gamma_A)=C^*_r(G_B)$ (or even a multiplier, if $f$ is bounded), not every such a function defines smooth (or even continuous) function on $B$ and only for such functions formula (\ref{kros-komut}) has immediate meaning.

{\em Comultiplication.} By computations similar to ones before  (\ref{def-Xr-kappa}) we obtain from (\ref{delta-Xr-1}):
\begin{equation*}
\delta(X^{A,r}_{\kropka{c}})(a_1 c_1, a_2  c_2)=\left((\tadc(a_1)\tadc(a_2) \kropka{c}) a_1 c_1\,, \, (\tadc(a_2)\kropka{c}) a_2 c_2\right).
\end{equation*}
and, for a basis $(X_\alpha)$ in $\gotc$: 
\begin{equation}
\delta(X^{A,r}_\alpha)(a_1 c_1, a_2 c_2)=\sum_\beta \tadc_{\beta \alpha}(a_2) X^{A,r}_\beta( a_1 c_1) + X^{A,r}_\alpha( a_2 c_2),
\end{equation}
or, as operators  on $\omh(\Gamma_A\times\Gamma_A)$  we can write (compare to (\ref{Delta0-gen}): \notka{Delta-gen-Gamma-A}
\begin{align}\label{Delta-gen-Gamma-A}
\Delta(A^A_\alpha)=I\mt A^A_\alpha+ \sum_\beta A^A_\beta \mt \tadc_{\beta \alpha}
\end{align}
But again, because of problems with domains, this expression is rather formal.

Finally, let us check that if we transfer formula (\ref{Delta-functions-B}) we get what we expect. The formula (\ref{wzor-embed}) gives diffeomorphism $\phi: A\rightarrow B'$:
\notka{dyfeo-A-Bprim}
\begin{equation}\label{dyfeo-A-Bprim}
\phi(a):= b_R(a)\,,\,\, \phi^{-1}(b')=a_R(b')\,,\,\,a\in A, b'\in B'
\end{equation}
For $f\in C^\infty(A)$ using (\ref{Delta-functions-B}) let us compute:
\begin{equation*}\begin{split}
[(\phi^*\mt \phi^*)^{-1}\Delta(\phi^*f)](a_1, a_2) &= \Delta(\phi^*f)](b_R(a_1), b_R(a_2))=(\phi^*f)(b_R[b_R(a_1)a_R(b_R(a_2))])=\\
&= (\phi^*f)(b_R(b_R(a_1)a_2))=(\phi^*f)(b_R(a_1a_2))=f(a_R(b_R(a_1a_2)))=\\
&= f(a_1 a_2)\end{split}
\end{equation*}

\subsection{Hopf $*$-algebra of $\kappa$-Poincar\'{e} from quantization of Poisson-Lie structure.}
On a  Hopf *-algebra level, relations for $\kappa$-Poincare Group were given in \cite{SZ-94} (and \cite{Kosinski-Maslanka}). The $*$-algebra  is generated by  {\em self-adjoint elements} 
$\{a_\alpha, L_{\alpha\beta}\,,\alpha,\beta=0,1\dots,n\}$ that satisfy the following commutation relations for $\eta_{\alpha\beta}:=diag(1,-1,\dots,-1)$ and some $h\in\R$ 
(to compare with \cite{Kosinski-Maslanka} substitute $\kappa:=h^{-1}$): \notka{hopf-rel-pedy}
%%%
\begin{align}\label{hopf-rel-pedy}
 [a_0,a_k] &= i h a_k\,,& \, [a_k,a_l] &=0 \,,& \,[L_{\alpha\beta}, L_{\gamma\delta}] &=0\,,& \,\eta=L^t\eta L
\end{align}\notka{hopf-rel-kros}
\begin{align}\label{hopf-rel-kros}
[a_0, L_{00}]& =i h ((L_{00})^2-1)\,,& \, [a_k, L_{00}]& = i h (L_{00}-1)L_{k0}\nonumber\\
[a_0, L_{0m}]& = i h L_{00} L_{0m}\,,& \, [a_k, L_{0m}]& = i h (L_{00}-1)L_{km}\nonumber\\
[a_0, L_{m0}]& = i h L_{00} L_{m0}\,,& \, [a_k, L_{m0}]& = i h (L_{k0} L_{m0}-\delta_{km}(L_{00}-1))\\
[a_0, L_{mn}]& =i h L_{m0} L_{0n}\,,& \, [a_k, L_{mn}]& = i h (L_{m0} L_{kn}-\delta_{km}L_{0n})\nonumber
\end{align}
together with coproduct $\Delta$:\notka{hopf-delta}
\begin{align}\label{hopf-delta}
\Delta(a_\alpha) &=a_{\alpha}\mt I + \sum_\beta L_{\alpha\beta}\mt a_{\beta}\,,& \Delta(L_{\alpha\beta}) &=\sum_\gamma L_{\alpha\gamma}\mt L_{\gamma\beta}
\end{align}
antypode $\antyp$  and counit $\counit$:\notka{hopf-antyp}\notka{hopf-counit}
\begin{align}
\antyp(a_\alpha) &=-\sum_\beta \antyp(L_{\alpha\beta}) a_\beta \,,& \antyp(L_{\alpha\beta}) &=(L^{-1})_{\alpha\beta}=(\eta L^t\eta)_{\alpha\beta}\label{hoph-antyp}\\
\counit(a_\alpha)&=0 \,,& \counit(L_{\alpha\beta}) &=\delta_{\alpha\beta}\label{hopf-counit}
\end{align}

\subsection{Comparison of formulae} Now we want to compare our generators and relations to formulae (\ref{hopf-rel-pedy} - \ref{hopf-delta}). 
We restrict computations only to commutation relations and comultiplication (and partially to antipode)  and are not going to investigate in details remaining parts of 
quantum group structure on $(C^*_r(\Gamma_A),\Delta)$  
(i.e. counit, antipode and Haar weight; this can be done with the use of expressions given in \cite{PS-DLG}).

Comparing  (\ref{kappa-e}) and  (\ref{Delta-gen-Gamma-A}) we see that, with respect to the comultiplication, our generators 
$(\widehat{Y}_\alpha)$ (where we put $\widehat{Y}_0:=\widehat{S}$) behave like $\antyp(a_\alpha)$. Such an identification would not be consistent with self-adjointness, however. 
Instead we use the unitary part of antipode (see e.g. \cite{SLW-MUQG}), which, for a quantum group defined by  DLG,  is implemented by the group inverse  \cite{PS-DLG}. 
As said above, we are not going to present all functional analytic details  but summarize  computations in the following:
\begin{lem} \label{lemat-Rgen}\notka{lemat-Rgen} Let $(\sM,\Delta)$ be a $*$-bialgebra generated by self-adjoint elements $(a_k, V_{kl}, V_{kl}^c)\,$, $k,l=1,\dots, n$ satisfying
$\displaystyle (1) \,\sum_m V_{km} V^c_{lm} = \sum_m V^c_{mk} V_{ml}=\delta_{kl} I\,$, $\displaystyle (2)\,\,\Delta(V_{kl}) =\sum_m V_{km}\mt V_{ml}\,$ and 
$\displaystyle (3)\,\,\Delta(a_k) = a_k\mt I +\sum_l V_{kl}\mt a_l.$ Then
\begin{itemize}
\item[a)] $\displaystyle \sum_m V^c_{lm} V_{km} = \sum_m V_{ml} V^c_{mk} =\delta_{kl} I$ and $\displaystyle \Delta(V_{kl}^c) =\sum_m V_{km}^c\mt V_{ml}^c.$
\item[b)] Elements $\,\displaystyle A_k := -\frac12 \sum_m \left( a_m V^c_{mk}+ V^c_{mk}a_m\right)$
are self-adjoint and satisfy
\begin{align}
\Delta(A_k) &= I\mt A_k +\sum_m A_m\mt V^c_{mk}
\end{align}
\item[c)] Moreover, if $\,\displaystyle  \sum_{mk}V^c_{mk} [a_m, V_{lk}] = \sum_{mk}[a_m, V_{lk}] V^c_{mk} \,$
then \notka{lemat-Rgen-eq1}
\begin{align}\label{lemat-Rgen-eq1}
a_k & = -\frac12 \sum_m\left( V_{km}A_m+A_m V_{km}\right)
\end{align}
and $\sM$ is generated by $(A_k,  V_{kl}, V_{kl}^c)$.
\end{itemize}
\end{lem}
\noindent {\em Proof:} a) The first equality is just $*$  applied to (1); for the second one, apply $\Delta$  to (1), use (2) and then (1); 
b) self-adjointness is evident and the formula for $\Delta(A_k)$ is a simple computation;\\
c) Clearly, equality (\ref{lemat-Rgen-eq1}) is sufficient for statement about generation. Let us rewrite the assumption in (c) as:
$$\sum_{mk}V^c_{mk} a_m V_{lk}- \sum_m \left(\sum_{k}V^c_{mk} V_{lk}\right) a_m = \sum_{m} a_m \left(\sum_k V_{lk}  V^c_{mk}\right) - \sum_{mk} V_{lk} a_m V^c_{mk}\,$$
by (1) and (a) expressions in brackets are $\delta_{lm} I$ and we obtain:\notka{lemat-Rgen-eq2}
\begin{equation}\label{lemat-Rgen-eq2}
2 a_l = \sum_{mk}\left( V^c_{mk} a_m V_{lk} + V_{lk} a_m V^c_{mk}\right)
\end{equation}
Let us compute using definition of $A_k$ and (1):
$$\sum_k V_{sk} A_k=-\frac12 \sum_{mk} V_{sk} a_m V^c_{mk} -\frac12 \sum_{mk} V_{sk} V^c_{mk}a_m=-\frac12 \sum_{mk} V_{sk} a_m V^c_{mk}-\frac12 a_s$$
In the similar way, using (a):
$$\sum_k  A_k V_{sk}=-\frac12 a_s - \frac12 \sum_{mk} V_{mk}^c a_m V^c_{sk}$$
Adding these two equalities and using (\ref{lemat-Rgen-eq2}) we get (\ref{lemat-Rgen-eq1}).
\dowl

The   version of the lemma above, with $a_k$ and $A_k$ interchanged,  is proven in the same way:\notka{lemat-Rgen1}
\begin{lem}\label{lemat-Rgen1}
Let $(\sM,\Delta)$ be a $*$-bialgebra generated by self-adjoint elements $(A_k, V_{kl}^c, V_{kl})\,$, $k,l=1,\dots, n$ satisfying
$\displaystyle (1) \,\sum_m V_{lm}^c V_{km} = \sum_m V_{ml} V_{mk}^c=\delta_{kl} I\,$, $\displaystyle (2)\,\Delta(V_{kl}^c) =\sum_m V_{km}^c\mt V_{ml}^c\,$ and 
$\displaystyle (3)\,\Delta(A_k) = I\mt A_k +\sum_l A_l \mt V_{lk}^c.$ Then
\begin{itemize}
\item[a)] $\,\displaystyle \sum_m V_{km} V^c_{lm} = \sum_m V^c_{mk} V_{ml}=\delta_{kl} I\,$ and  $\displaystyle \,\Delta(V_{kl}) =\sum_m V_{km}\mt V_{ml}\,$;
\item[b)] Elements $\,\displaystyle a_k  = -\frac12 \sum_m\left( V_{km}A_m+A_m V_{km}\right)$
are self-adjoint and satisfy
\begin{align}
\Delta(a_k) = a_k\mt I +\sum_l V_{kl}\mt a_l.
\end{align}
\item[c)] Moreover, if $\,\displaystyle  \sum_{mk} [V^c_{ks}, A_m] V_{km} = \sum_{mk}V_{km} [V^c_{ks}, A_m]\,$
then \notka{lemat-Rgen1-eq1}
\begin{align}\label{lemat-Rgen1-eq1}
A_k &= -\frac12 \sum_m \left( a_m V^c_{mk}+ V^c_{mk}a_m\right)
\end{align}
and $\sM$ is generated by $(a_k,  V_{kl}^c, V_{kl})$.
\end{itemize}
\end{lem}\dowl

Looking at the formulae (\ref{Delta-na-generatorach}) and using the fact that $W$ given by (\ref{kad-zUd}) is a representation of the group $A$ we use the lemma \ref{lemat-Rgen1} with 
$V^c=W$, where matrix elements of $W$ are expressed by $(\Lambda,w,u,\alpha)$ as in (\ref{kad-zUd-ar}) (using identification of $A$  with $B'$ as in (\ref{dyfeo-A-Bprim})).
So we define:\notka{def-L}
\begin{equation}\label{def-L}
L:=\left(\begin{array}{ccc}
\frac{1}{\alpha} & - \frac{w^t}{\alpha} \\ 
\frac{u}{|\alpha|}  & sgn(\alpha)( \Lambda-\frac{u w^t}{\alpha}) \end{array}\right)
\end{equation}
and \notka{def-a}
\begin{equation}\label{def-a}
\begin{split}
a_\alpha &:= -\frac12 \sum_\beta\left( L_{\alpha\beta}\widehat{Y}_\beta+\widehat{Y}_\beta L_{\alpha\beta}\right)=
-\sum_\beta L_{\alpha\beta}\widehat{Y}_\beta -\frac12 \sum_\beta [\widehat{Y}_\beta , L_{\alpha\beta}]=\\
&=: \tilde{a}_\alpha -\frac12 \sum_\beta [\widehat{Y}_\beta , L_{\alpha\beta}]
\end{split}
\end{equation}
In the formula above and in what follows we treat $\widehat{Y}_\alpha$ (recall that $\widehat{Y}_0 :=\widehat{S}$) 
and matrix elements of $L$ {\em as operators on $\omh(\Gamma_A)$ (or $\sA(\Gamma_A)$)}.

\begin{re}
Notice that matrix elements of $L$, despite being non-continuous functions on $B$, are affiliated to $C^*(G_B)$ (wich is equal to  $C^*_r(G_B)$)  
because they are smooth functions on $A$ and $C^*(G_B)=C^*_r(\Gamma_A)$.
\end{re}

\noindent
Since for $f\in C^\infty(A)$  commutators $[\widehat{Y}_\mu , f]\in C^\infty(A)$  and functions on $A$ commute, we have
$$[a_\mu , L_{\beta\gamma}]=[\tilde{a}_\mu ,  L_{\beta\gamma}]$$

We will need commutators of $\widehat{Y}_\beta$ with $\frac{1}{\alpha}, sgn(\alpha)$ and $\frac{1}{|\alpha|}$. Since $\alpha\neq 0$ on $A$ and $\widehat{Y}_\beta$ are differential operators, 
it is clear that \notka{rel-SY-sgn}
\begin{equation}\label{rel-SY-sgn}
[\widehat{Y}_\beta, sgn(\alpha)]=0.
\end{equation}
By this equality we have 
$\,\displaystyle \left[\widehat{Y}_\beta, \frac{1}{|\alpha|}\,\right]=$ $ \displaystyle \left[\widehat{Y}_\beta , \frac{sgn(\alpha)}{\alpha}\,\right]=$ 
$ \displaystyle sgn(\alpha) \left[\widehat{Y}_\beta, \frac{1}{\alpha}\,\right]$
and, for the last commutator, 
$$0=\left[\widehat{Y}_\beta, \frac{1}{\alpha} \alpha\right]=\left[\widehat{Y}_\beta, \frac{1}{\alpha}\right] \alpha+ \frac{1}{\alpha} \left[\widehat{Y}_\beta, \alpha\right]
\quad \Rightarrow \quad \left[\widehat{Y}_\beta, \frac{1}{\alpha}\right]=- \frac{1}{\alpha} \left[\widehat{Y}_\beta, \alpha\right] \frac{1}{\alpha}. $$
Now, using (\ref{relkom-S-1}), (\ref{relkom-Y-1}) we obtain: 
\begin{align*}
[\widehat{S}, \frac{1}{\alpha}\,] &=  \iota(1-\frac{1}{\alpha^2})\,,&\,\,[\widehat{Y}_m, \frac{1}{\alpha}\,] &= \iota\frac{\alpha-1}{\alpha^2} w_m\,.
\end{align*}
%%%%%%%%%
Having  these relations together with  (\ref{relkom-S-1}), (\ref{relkom-Y-1}) by direct computation on obtains commutators of $\widehat{Y}_\gamma$ with matrix elements of $L$:
\notka{relkom-SY-L}
\begin{equation}\label{relkom-SY-L}
\left[\widehat{Y}_\gamma, L_{\beta\mu}\right]=
\iota\left( \delta_{\gamma \mu}(\delta_{\beta 0}-L_{\beta 0})- sgn(\gamma) L_{\beta \gamma}(L_{0\mu}-\delta_{0 \mu})\right),
\end{equation}
where we use $sgn(\gamma):=\left\{\begin{array}{lcr} 1&\,{\rm for}\, & \gamma=0\\
                                                      -1 & \, {\rm for}\, & \gamma>0\end{array}\right..$
Now we get
\begin{equation}
\begin{split}
\left[a_\rho, L_{\beta \mu}\right] &= \left[\tilde{a}_\rho, L_{\beta \mu}\right]=- \sum_\gamma\left[L_{\rho \gamma}\widehat{Y}_\gamma, L_{\beta \mu}\right]=
 - \sum_\gamma L_{\rho \gamma} \left[\widehat{Y}_\gamma, L_{\beta \mu}\right]= \\
&= \iota\left( L_{\rho \mu} (L_{\beta 0}-\delta_{\beta 0}) +sgn(\rho) \delta_{\rho\beta}(L_{0 \mu} - \delta_{0 \mu})\right).
\end{split}\end{equation}
and these are relations (\ref{hopf-rel-kros}) (for $h=1$).

Finally, lets us verify (\ref{hopf-rel-pedy}). By (\ref{relkom-SY-L}) and (\ref{relkom-S-Y}):
\begin{equation*}
\sum_\gamma[\widehat{Y}_\gamma, L_{\beta \gamma}] = \iota n (\delta_{\beta 0}- L_{\beta 0})
\end{equation*}
\begin{equation*}
\begin{split}
[\tilde{a}_\mu, \tilde{a}_\nu] &= \sum_{\gamma \delta} \left[ L_{\mu \gamma} \widehat{Y}_\gamma, L_{\nu \delta} \widehat{Y}_\delta \right]= 
\sum_{\gamma \delta} L_{\mu \gamma} L_{\nu \delta}  \left[  \widehat{Y}_\gamma, \widehat{Y}_\delta \right] + 
 L_{\mu \gamma} \left[  \widehat{Y}_\gamma, L_{\nu \delta} \right] \widehat{Y}_\delta 
- L_{\nu \delta}  \left[ \widehat{Y}_\delta , L_{\mu \gamma} \right]  \widehat{Y}_\gamma=\\
&= \iota(\delta_{0\mu}\tilde{a}_\nu - \delta_{0\nu}\tilde{a}_\mu).
\end{split}
\end{equation*}
and
\begin{equation}
\begin{split}
[a_\mu, a_\nu] &= \left[\tilde{a}_\mu- \frac{\iota n}{2} (\delta_{\mu 0}- L_{\mu 0}), \tilde{a}_\nu- \frac{\iota n}{2} (\delta_{\nu 0}- L_{\nu 0})\right]= \\
 &= \left[\tilde{a}_\mu , \tilde{a}_\nu \right]+ \frac{\iota n}{2}\left( \left[\tilde{a}_\mu ,  L_{\nu 0})\right]- \left[\tilde{a}_\nu, L_{\mu 0})\right]\right)
= \iota(\delta_{0\mu}\tilde{a}_\nu - \delta_{0\nu}\tilde{a}_\mu)+ \frac{\iota n}{2} \iota \left( \delta_{\mu 0} L_{\nu 0} - \delta_{\nu 0} L_{\mu 0} \right) =\\
&= \iota \left( \delta_{\mu 0} (\tilde{a}_\nu + \frac{\iota n}{2} L_{\nu 0}) - \delta_{\nu 0} (\tilde{a}_\mu + \frac{\iota n}{2} L_{\mu 0})\right)= \\
 &= \iota \left( \delta_{\mu 0} (\tilde{a}_\nu + \frac{\iota n}{2} (L_{\nu 0} -\delta_{\nu 0})) - \delta_{\nu 0} (\tilde{a}_\mu + \frac{\iota n}{2} (L_{\mu 0}-\delta_{\mu 0}) \right)=\\
&= \iota \left( \delta_{\mu 0} a_{\nu} - \delta_{\nu 0} a_\mu \right)
\end{split}
\end{equation}
as in (\ref{hopf-rel-pedy}).

%%%%%%%%%%%%%%%%%%%%%%%%%%%%%%%%%%%%%%%%%%%%%%%%%%%%%%%%%%%%%%%%%%%%%%%%%%%%%%%%%%%%%%%%%%%%%%%%%%%%%%%%%%%%%%%%%%
%%%%%%%%%%%%%%%%%%%%%%%%%%%%%%%%%%%%%%%%%%%%%%%%%%%%%%%%%%%%%%%%%%%%%%%%%%%%%%%%%%%%%%%%%%%%%%%%%%%%%%%%%%%%%5

\section{(Some remarks on) Quantum $\kappa$-Minkowski space.}
In this last section we are going to make few remarks on ``quantum $\kappa$-Minkowski space''. 
Under this name is usually understood  ``$*$-algebra generated by self-adjoint elements $(x_0, x_k)$ satisfying relations''  (\ref{hopf-rel-pedy}):
\begin{align*}%\label{hopf-rel-pedy}
 [x_0,x_k] &= i h x_k\,,& \, [x_k,x_l] &=0 \quad,\quad k=1,\dots,n
\end{align*}
As a Hopf algebra it was introduced already in \cite{Maj-Rueg}. This is considered either on a pure algebra level or as operators on Hilbert space \cite{Agostini} 
or within a star-product formulation of simplified two-dimensional version in \cite{Sitarz}.   But if one wants to consider it on a $C^*$-algebra level, it is clear 
(in fact since   \cite{PS-triple}, \cite{SZ-PP} and certainly since \cite{VV}) that  the ``unique candidate'' for this name is $C^*(C)$  with the group $C$ defined in (\ref{def-C}) i.e. 
$C$ is the $AN$ group from the Iwasawa decomposition $SO_0(1,n)=SO(n) A N$ (this group appears in \cite{Agostini} under the name  ``$\kappa$-Minkowski Group''). 
But to call a space   the ``Minkowski space'' it should carry an action of the Poincar\'{e} Group, so in our case we have to show the action of our quantum 
group $(C^*(\Gamma_A), \Delta)$ on $C^*(C)$.
Let us now show, how the groupoid framework gives natural candidate for such an action and postpone analytic details to future publication.
\begin{re}
A ``quantum space'' $C^*(C)$ has a family of classical points $\{p_\lambda:\lambda\in\R\}$ i.e. characters. They are given by 1-dimensional unitary representations of $C$
$p_\lambda(s,y):=s^{\iota\lambda}$. These are representations induced from whole line of $0$-dimensional symplectic leaves on which the Poisson bivector  for the related Poisson Minkowski 
space (affine) described e.g.  in \cite{PS-poisson} vanishes. 
\end{re}
Let us begin with simpler situation of global decomposition and 
let $(G;B,C)$ be a double group with $\delta_0: G_B\rel G_B\times G_B$ defined by (\ref{def-delta0}).  
Consider relations $\delta_L: C\rel G_B\times C$ and $\delta_R: C\rel C\times G_B$ defined by:\notka{def-deltaLR}
\begin{align} \label{def-deltaLR} 
\delta_L &:= \{(g,c_R(g); c_L(g))\,:\,g\in G\}\quad,\quad 
\delta_R:=\{(c_L(g), g; c_R(g))\,:\,g\in G\}
\end{align}
The following lemma can be proven by direct computation:
\begin{lem} $\delta_L$ and $\delta_R$ are morphisms of groupoids that satisfy:
\begin{align}
(\delta_0\times id)\delta_L &= (id\times \delta_L)\delta_L\quad, & \quad ( id \times \delta_0)\delta_R &= (\delta_R\times id )\delta_R
\end{align}
Let $\sigma:C\times G_B \ni (c,g)\mapsto (g,c)\in G_B\times \C$ be the flip and $R_B, R_C$ denote the inverse in a group $G$ and  its restriction to $C$ respectively. 
Then the following equality holds:
\begin{align} 
\delta_L &= \sigma(R_C\times R_B)\delta_R R_C
\end{align}
 \end{lem}
\dowl

If $(G;B,C)$ is a DLG, relations $\delta_L,\delta_R$ are morphisms of differential groupoids and lifting them  one  obtains $\Delta_L,\Delta_R$ -- morphisms of corresponding $C^*$-algebras,
e.g. $\Delta_L\in Mor(C^*(C), C^*(G_B\times C))$. One may expect that in  ``nice'' cases $C^*(G_B\times C)=C^*(G_B)\otimes C^*(C)$ and reduced/universal algebras problem can be also handled. 
This way one gets actions (left or right) of $(C^*(G_B),\Delta_0)$ on  $C^*(C)$.  Let us show that $\Delta_L$ should be considered as ``quantization'' of the canonical affine action
of the semidirect product $\gotb^0\rtimes B$ on $\gotb^0$:\notka{affine-action}
%%%%%%%%%%%%%%%%%%%%
\begin{equation}\label{affine-action}
(\gotb^0\rtimes B)\times \gotb^0\ni (\varphi,b;\psi)\mapsto \varphi+\kad(b)\psi\in \gotb^0
\end{equation}

Let $\Gamma_1\rightrightarrows E_1$ and $\Gamma_2\rightrightarrows E_2$ be differential groupoids. For  a (differentiable) relation $h:\Gamma_1\rel\Gamma_2$ 
by $T^*h: T^*\Gamma_1\rel T^*\Gamma_2$ we  denote the relation ({\em cotangent lift of $h$}) defined by \notka{def-Twgiazdkah}
\begin{equation} \label{def-Tgwiazdkah}
(\varphi_2,\varphi_1)\in T^*h \iff   \forall (v_2,v_1) \in Th \quad  <\varphi_2,v_2 > = <\varphi_1,v_1 >.
\end{equation}
%%%%%%%%%%%%%%%%%
If $h$ is a morphism of differential groupoids then $T^*h$ is a morphism of {\em symplectic groupoids} and its base map is a Poisson map $(TE_2)^0\rightarrow (TE_1)^0$ \cite{SZ2}.

By (\ref{def-deltaLR}), (\ref{def-Tgwiazdkah}) and (\ref{kontragrad}) for $\Phi$ -- the base map of $T^*\delta_L$, we have $\Phi: (TB)^0\times \gotc^*\rightarrow \gotc^*$ and 
for $(\phi, \psi)\in (T_bB)^0\times \gotc^*$, $\kropka{c}\in\gotc\subset\gotg$ and identifying $\gotc^*$ with $\gotb^0\subset\gotg^*$:
%%%%%%%%%%%%%
\begin{equation*} 
\begin{split}
<\Phi(\phi, \psi), \kropka{c}> &=  <(\phi, \psi), (\kropka{c} b, c_R(\kropka{c}b)>=<\phi,\kropka{c}b > + < \psi , \adc(b^{-1})(\kropka{c})>= \\
 &= < \phi,\kropka{c}b > + < \kad(b)(\psi) ,\kropka{c})>= <\varphi+\kad(b)(\psi),\kropka{c}>,
\end{split}
\end{equation*}
where, using right trivialization, we represent $\phi\in (T_bB)^0$ by $(\varphi, b)\in \gotb^0\times B$. This way we get
$\quad\displaystyle \Phi(\varphi,b; \psi)=\varphi+\kad(b)(\psi)\quad$ exactly as in (\ref{affine-action}).

In the situation of $\kappa$-Poincar\'{e}, relations defined as $\delta_L,\delta_R$ in (\ref{def-deltaLR}) e.g. 
$$\tilde{\delta}_L: C\rel \Gamma_A \times C \quad,\quad \tilde{\delta}_L:=\{(g,\tilde{c}_R(g);\tilde{c}_L(g)) :g\in\Gamma_A\}$$
{\em are  not morphisms} of differential groupoids. It can be directly verified that $\quad\displaystyle m'(\tilde{\delta}_L \times \tilde{\delta}_L)\nsubseteq \tilde{\delta}_L m\,\,$, where 
$m$ and $m'$ denote multiplication relations in $C$ and $\Gamma_A \times C$, respectively. Or one may observe that $T^*\tilde{\delta}_L$ restricted to sets of units gives the action of 
Poisson-Poincar\'{e} group $\gota^0\rtimes A$ on $\gotc^*\simeq \gota^0$  as 
in (\ref{affine-action}). It is known that this is a Poisson action which is {\em non complete}, so it cannot be the base map of a morphism of symplectic groupoids \cite{SZ2}. 
Essentially we face the similar problem as in the very beginning: some operators that ``should be'' self-adjoint  are not essentially self-adjoint on their ``natural'' domains and we can 
try to overcome it in the similar way -- passing from $\Gamma_A$ to $G_B$.

It is easier to work with $\tilde{\delta}_R:C\rel C\times \Gamma_A$ given by 
$\displaystyle \,\tilde{\delta}_R:=\{(\tilde{c}_L(g), g; \tilde{c}_R(g))\,:\,g\in \Gamma_A\}\,,\,$ or after passing to  $\Gamma_{B'}$  (and using the same symbol $\tilde{\delta}_R$):
\notka{def-tilde-deltaR}
\begin{equation}\label{def-tilde-deltaR}
\tilde{\delta}_R=\{(\tilde{c}_L(b_L(g))^{-1}\tilde{c}_L(g), g; c_R(g)) : g\in \Gamma_{B'}\}\subset (C\times \Gamma_{B'}) \times C
\end{equation}
Let us define \notka{def-TC}
\begin{align}\label{def-TC}
T^C:=\{(\tilde{c}_L(b)^{-1},b):b\in B'\}=(c_R\times id)T\subset C\times G_B\,, \quad T^C_{12}:= T^C\times B\subset C\times G_B\times G_B
\end{align}
The following lemma may be compared to Prop. \ref{prop-twist}.\notka{lemma-TC}
\begin{lem}\label{lemma-TC}
\begin{enumerate}
\item $T^C$ is a section of left and right projections over $\{e\}\times B'\subset C\times G_B$ ($e$ is the neutral element in $C$) and a bisection of $C\times \Gamma_{B'}$.
\item $(id\times\delta_0) T^C$ is a section of left and right projections (in $C\times G_B\times G_B$) over the set $\{e\}\times \delta_0(B')=\{(e,b_2,b_3) : b_2 b_3\in B'\}$;
\item $(\delta_R\times id) T^C$ is a section of left and right projections  over the set $\{e\}\times B \times B'$;
\item  $T^C_{12}(\delta_R\times id) T^C= T_{23}( id \times \delta_0) T^C$ (equality of sets in $C\times G_B\times G_B$), 
moreover this set is a section of the right projection over $\{e\}\times (\delta_0(B')\cap(B\times B'))$ and  the left projection over
$\{e\}\times B'\times B'$.
\end{enumerate}
\end{lem}

{\em Proof:} The first statement is clear from the definition of $T^C$.
By a straightforward computation:
\begin{equation}\begin{split}
(id\times\delta_0) T^C & =\{(\tilde{c}_L(b_1 b_2)^{-1}, b_1, b_2): b_1 b_2\in B'\}\subset C\times G_B\times G_B\\
(\delta_R\times id) T^C &= \{(c_L(g), g, b): b\in B', c_R(g)=\tilde{c}_L(b)^{-1}\}=\\
 &= \{(c_L(b_2 \tilde{c}_L(b_3)^{-1}), b_2 \tilde{c}_L(b_3)^{-1}, b_3): b_2\in B, b_3\in B'\}\subset C\times G_B\times G_B,
\end{split}
\end{equation}
 the second and third statement follow easily from these expressions.
Now, using these expressions, it is easy to compute $T^C_{12}(\delta_R\times id) T^C$ and get:
$$T^C_{12}(\delta_R\times id) T^C= \{\tilde{c}_L(b_2)^{-1} c_L(b_2 \tilde{c}_L(b_3)^{-1}), b_2 \tilde{c}_L(b_3)^{-1}, b_3) : b_2,b_3\in B'\},$$
and the same result for $T_{23}( id \times \delta_0) T^C$. The statement about left projection is clear, the only relation we need to check is that for $b_2,b_3\in B'$ the product 
$b_R(b_2 \tilde{c}_L(b_3)^{-1}) b_3$ is again in $B'$. Let $b_2=c_2 a_2$ and $b_3=c_3 a_3 $. Then
$$b_R(b_2 \tilde{c}_L(b_3)^{-1}) b_3=b_R(c_2 a_2 c_3^{-1})b_3= b_R(c_2 a_2 c_3^{-1}b_3)= b_R(c_2 a_2 a_3)=b_R(a_2 a_3)\in B'$$
\dowl

\noindent
$T^C$ is used to twist $\delta_R$ to $\tilde{\delta}_R$:
\begin{lem} Let $\delta_R:C\rel C\times G_B$ be the morphism defined in (\ref{def-deltaLR})  and $\tilde{\delta}_R$ be as in (\ref{def-tilde-deltaR}). 
They are related by: $\displaystyle \tilde{\delta}_R=Ad_{T^C}\cdot \delta_R$.
\end{lem}
{\em Proof:} Recall that $Ad_{T^C}:C\times G_B\rel C\times G_B$ is defined by (compare (\ref{def-AdT})):
$$(c_1,g_2;c_3,g_4)\in Ad_{T^C}\iff \exists t_1,t_2\in T^C : (c_1,g_2) = t_1 (c_3,g_4) (s_C\times s_B) t_2$$
The multiplication above is in the groupoid $C\times G_B$ and $s_C, s_B$ are groupoid inverses (i.e. $s_C$ stands for the inverse in the group $C$). Let us  compute:
$$(c_1,g_2;c_3,g_4)\in Ad_{T^C}\iff \exists b_5,b_6\in B' : (c_1,g_2)= (\tilde{c}_L(b_5)^{-1}, b_5) (c_3, g_4) (\tilde{c}_L(b_6), b_6),$$
i.e. $b_5=b_L(g_4)\in B'\,,\, b_6=b_R(g_4)\in B'$ therefore  $g_2=g_4\in \Gamma_{B'}$ and $c_1=\tilde{c}_L(b_5)^{-1} c_3 \tilde{c}_L(b_6);$
so this relation is in fact bijection
$$ Ad_{T^C}: C\times\G_{B'}\ni (c,g)\mapsto (\tilde{c}_L(b_L(g))^{-1} c_3 \tilde{c}_L(b_R(g)),g)\in C\times\G_{B'}$$
defined by the bisection $T^C$ of $C\times\G_{B'}$.

\noindent
Now we have:
$$(c_1,g_2,c_3)\in Ad_{T^C}\delta_R \iff \exists g\in G : (c_1,g_2; c_L(g), g)\in Ad_{T^C}\,, \,c_R(g)=c_3$$
therefore $g\in \G_{B'}$ and $c_1=\tilde{c}_L(b_L(g))^{-1} c_L(g) \tilde{c}_L(b_R(g))=\tilde{c}_L(b_L(g))^{-1} \tilde{c}_L(g)$ and
$$Ad_{T^C}\delta_R=\{(\tilde{c}_L(b_L(g))^{-1} \tilde{c}_L(g), g, c_R(g)) : g\in \Gamma_{B'}\}$$
exactly as $\tilde{\delta}_R$ in (\ref{def-tilde-deltaR}).\dowl\\

$T^C$, as a bisection of $C\times \Gamma_{B'}$, defines unitary multiplier $\widehat{T^C}$ of $C^*_r(C\times\Gamma_{B'})=C^*_r(C)\otimes C^*_r(\Gamma_{B'})=
C^*(C)\otimes C^*_r(G_B)=C^*(C)\otimes C^*(G_B)$, so having $\Delta_R$ -- action of the quantum group $(C^*(G_B),\Delta_0)$ on ``quantum space'' represented by 
$C^*(C)$ (lifted from $\delta_R$)  we can, due to properties of $T^C$ described in Lemma \ref{lemma-TC}, {\em define the action} of our quantum group 
$(C^*(\Gamma_A),\Delta)$ on the same ``space'' $C^*(C)$  by  $\widetilde{\Delta}_R(a):=\widehat{T^C}\Delta_R(a)\widehat{T^C}^{-1}$. 
Whether this construction gives {\em continuous}  action of  $\kappa$-Poincar\'{e} on $C^*(C)$ still needs to be verified.

\section{Appendix}

Here we collect some formulae proven in \cite{PS-DLG} and used in this paper.
$(G;B,C)$ is a double Lie group,   $\gotg, \gotb,\gotc$  are corresponding Lie algebras and $\gotg=\gotb\oplus \gotc$ 
(direct sum of vector spaces).
Let  $\rzutB, \rzutC $ 
be projections in $\gotg$ corresponding to the decomposition $\gotg=\gotb\oplus \gotc$.
Let us define:
\begin{equation}\label{adb-adc}\notka{adb-adc}
\adb(g):=\rzutB Ad(g)|_{\gotb}\,\,,\,\,\adc(g):=\rzutC  Ad(g)|_{\gotc}
\end{equation}
Clearly $\adb$ and $\adc$ are representations when restricted to $B$ or $C$.

\noindent{\em Modular functions.}
Let us define:
\notka{app-def-modular-functions}\begin{equation}\label{app-def-modular-functions}
j_B(g):=|\det(\adb(g))|\,,\,\,j_C(g):=|\det(\adc(g))|
\end{equation}

\noindent {\em The choice of $\om_0$.}
Choose a real half-density $\mu_0\neq 0$ on $T_eC$ and define left-invariant half-density on $G_B$ by
\begin{equation}\label{def-lo}\notka{def-lo}
\lambda_0(g)(v):=\mu_0(g^{-1} v)\,,\,v\in \lma T^l_g G_B.
\end{equation}
The corresponding right-invariant  half-density is given by:
\begin{equation}\label{def-ro}\notka{def-ro}
\rho_0(g)(w):=j_C(b_L(g))^{-1/2}\mu_0(w g^{-1})\,,\,w\in\lma T^r_g G_B.
\end{equation}
{\em Multiplication and comultiplication in $\sA(G_B)$}
After the choice of $\om_0$ as above,  the multiplication in $\sA(G_B)$ reads:
$(f_1\om_0)(f_2\om_0)=:(f_1*f_2)\om_0$ and
\notka{mult-ga}
\begin{equation}\label{mult-ga}
\begin{split}
(f_1*f_2)(g) & =\int_C d_lc\, f_1(b_L(g) c)f_2(c_L(b_L(g) c)^{-1}g)=\\
&=\int_C d_rc\, j_C(b_L(c b_R(g)))^{-1} f_1(g  c_R(c b_R(g))^{-1}) f_2(c b_R(g)),
\end{split}
\end{equation}
where $d_l c$ and $d_r c$ are left and right Haar measures on $C$ defined by $\mu_0$.

The $||\cdot||_l$ defined by this $\om_0$ is given by:
\notka{l-norm}\begin{equation}\label{l-norm}
||f||_l=\sup_{b\in B}\int_C d_l c\, |f(b c)|
\end{equation}
Let $\delta_0:=m_C^T: G_B\rel G_B\times G_B$. The formula for $\hat{\delta}_0$ reads
$$\hat{\delta}_0(f\om_0)(F(\om_0\otimes\om_0))=:(\hat{\delta}_0(f) F)(\om_0\otimes \om_0)$$
\notka{app-hat-delta0-form}\begin{equation}\label{app-hat-delta0-form}
(\hat{\delta}_0(f) F)(b_1 c_1, b_2 c_2)=
\int_C d_lc\, j_C(c_L(b_2 c))^{-1/2}f(b_1 b_2 c)F(c_L(b_1 b_2 c)^{-1} b_1 c_1, b_R(b_2 c) c^{-1} c_2)
\end{equation}

\noindent{\em Action of bisections on bidensities.} Let us write $\om=f \om_0\,,\,f\in \sD(G_B)$. The action of a bisection $B c_0$ on $\sA(G_B)$ is given by:
\notka{app-Bc0-action}
\begin{align}\label{app-Bc0-action}
 (B c_0)(f\om_0) & =:(B c_0 f)\, \om _0\,, &\quad & (B c_0 f)(g)=f(B(c_0)^{-1} g) j_C(c_0)^{-1/2}
\end{align}

\noindent
{\bf Notation for orthogonal Lie algebras.}
Let $(V,\eta)$ be a  real, finite dimensional vector space with a bilinear, symmetric  and non degenerate  form $\eta$;
by $\eta$ we  denote also the isomorphism $V\rightarrow V^*$ defined  by $< \eta(x), y >:=\eta(x,y)$.
A basis $(v_\alpha)$ of  $V$ is called {\em orthonormal} if 
$\displaystyle \eta(v_\alpha, v_\beta)=\eta(v_\alpha,v_\alpha)\delta_{\alpha\beta}\,,\,\,\,|\eta(v_\alpha,v_\alpha)|=1$.
For a subset $S\subset V$ the symbol $S^\perp$  is used for {\em the orthogonal complement} of $S$, 
the symbol $S^0\subset V^*$ stands for {\em the annihilator} of $S$.
%%%%%%%%%%%%%%
Let us define operators in $End(V)$:
\begin{equation}\label{lambdas}
\Mlambda_{xy}:=x\mt\eta(y)-y\mt\eta(x)\quad, \quad x,y\in V.\end{equation} 
For a basis $(v_\alpha)$ in $V$ we  write $\Mlambda_{\alpha\beta}$ instead of  $\Mlambda_{v_\alpha,v_\beta}$. 
Operators $\Mlambda_{xy}$ satisfy:
\begin{equation}\label{lambda-komut} \notka{lambda-komut} 
\,[\Mlambda_{xy},\Mlambda_{zt}]=\eta(x,t)\Mlambda_{yz}+\eta(y,z)\Mlambda_{xt}-\eta(x,z)\Mlambda_{yt}-\eta(y,t)\Mlambda_{xz}
\end{equation}
and $so(\eta)=span\{\Mlambda_{xy}:x,y\in V\}$.
Notice   that for  $g\in O(\eta)$ we have  $\ad(g) (\Mlambda_{xy})=\Mlambda_{gx,gy}$. 
We will use a bilinear, non degenerate form $k:so(\eta)\times so(\eta)\rightarrow \R$ defined by: \notka{app-def-k}
\begin{equation}
\label{app-def-k}
k(\Mlambda_{xy},\Mlambda_{zt}):=\eta(x,t)\eta(y,z)-\eta(x,z)\eta(y,t)
\end{equation}
It is easy to see that for $g\in O(\eta)$: $\ad(g)\in O(k)$ i.e.
$$k(g\Mlambda_{xy}g^{-1},g\Mlambda_{zt}g^{-1})=k(\Mlambda_{xy},\Mlambda_{zt})\,,\,g\in O(\eta)$$
(of course $k$ is proportional to the Killing form on $so(\eta)$). By $\kad$ we denote the coadjoint representation
of $O(\eta)$ on $so(\eta)^*$: $\kad(g):=\ad(g^{-1})^*$.
If $k$ is the isomorphism $so(\eta)\rightarrow so(\eta)^*$ defined by the form $k$ then 
$$\kad(g) k(X)=k(\ad(g) X) \,,\quad X\in so(\eta), g\in O(\eta)$$
Let us also define  a bilinear form $\tilde{k}$ on $so(\eta)^*$ by: 
\begin{equation}\label{tildek}\notka{tildek}
\tilde{k}(\varphi,\psi):=k(k^{-1}(\varphi),k^{-1}(\psi))\,,\,\varphi,\psi\in so(\eta)^*
\end{equation}
so $\tilde{k}(\varphi,\psi)=<\varphi, k^{-1}(\psi)>$;  
again it is clear that if $g\in O(\eta)$ then $\kad(g)\in O(\tilde{k})$, and
$$\tilde{k}(\kad(g) k(X),k(Y))= k(\ad(g) X, Y)\,\,,\,X,Y\in so(\eta).$$

\noindent
{\bf Adjoint, coadjoint representations and Hopf algebra structure.} For a Lie group $G$ we identify the group $T^*G$, via right translations, with the semidirect product
$\gotg^*\rtimes G$ (with coadjoint representation):\notka{app-semi-direct}
\begin{equation}
\label{app-semi-direct}
(\varphi,g)(\psi, h):=(\varphi+\kad(g) \psi, g h)\,,\quad \varphi, \psi\in \gotg^*\,,\,g,h\in G,
\end{equation}
If  $B\subset G$ is a subgroup with a Lie algebra $\gotb\subset \gotg$ then $\gotb^0\times B$ is  a subgroup of $\gotg^*\rtimes G$.

If $\gotc\subset \gotg$ is any complementary subspace to $\gotb$ i.e. $\gotg=\gotb\oplus\gotc$, then $Ad^{\gotc}(b):=P_{\gotc}Ad(b)|_{\gotc}$ is a representation of $B$ on $\gotc$.
The spaces $\gotc$ and $\gotb^0$ are dual to each other and the representation $Ad^{\gotc}$ is contragradient to $\kad|_B$, i.e. for $\varphi\in \gotb^0, \kropka{c}\in\gotc$ and $b\in B$:
\notka{kontragrad}
%%%%%%%%%%%%%%%%
\begin{equation}\label{kontragrad}
\begin{split}
<\kad(b)\varphi, \kropka{c}> &= <\varphi, Ad(b^{-1})\kropka{c}>=<\varphi, P_{\gotc} Ad(b^{-1})\kropka{c}>=<\varphi, Ad^{\gotc}(b^{-1})\kropka{c}>=\\
&= <(Ad^{\gotc}(b^{-1}))^*\varphi, \kropka{c}>
\end{split}
\end{equation}
Let $(\rho_k)$ be a basis in $\gotb^0$,  $(\kropka{c}_k)$ dual basis in $\gotc$ and   $\kad_{lk}, Ad^{\gotc}_{lk} : B\rightarrow\R$  matrix elements of  $\displaystyle \kad(b)$ and $Ad^{\gotc}(b)$ 
in corresponding bases:
\begin{equation*} 
<\rho_l, \kropka{c}_m>=\delta_{lm} \quad,\quad\quad \kad(b)\rho_k=\sum_l \kad_{lk}(b)\rho_l\quad,\quad\quad Ad^{\gotc}(b)\kropka{c}_k=\sum_l Ad^{\gotc}_{lk}(b) \kropka{c}_l
\end{equation*}
Clearly, the equality  (\ref{kontragrad}) implies $\,\displaystyle \kad_{lk}(b^{-1})=Ad^{\gotc}_{kl}(b)$, i.e.
$$\sum Ad^{\gotc}_{km}(b)\kad_{kl}(b)=\sum \kad_{lk}(b)Ad^{\gotc}_{mk}(b)=\delta_{lm}.$$
Let us  use the same symbols for extensions of functions $\kad_{lk}, Ad^{\gotc}_{lk},  \kropka{c}_k$ to  $\gotb^0\rtimes B$ i.e.
$$\kad_{kl}(\varphi,b):=\kad_{kl}(b)\,,\,\,Ad^{\gotc}_{kl}(\varphi,b):=Ad^{\gotc}_{kl}(b)\,,\,\,\kropka{c}_k(\varphi,b):=<\varphi, \kropka{c}_k>$$
%%%%%%%%%%
It is straightforward to compute action of comultiplication, counit and antipode, defined by the group $\gotb^0\rtimes B$, on functions $\kropka{c}_k, \kad_{kl}, Ad^{\gotc}_{kl}$:\notka{hopf}
\begin{align*}
\Delta(\kropka{c}_k)  &= \kropka{c}_k\mt I + \sum_m \kad_{km}\mt \kropka{c}_m\,, & \counit(\kropka{c}_k) &= 0\,,   & \antyp(\kropka{c}_k) = -\sum_m\antyp(\kad_{km})\,\kropka{c}_m &= -\sum_m Ad^{\gotc}_{mk} \,\kropka{c}_m\,;
\end{align*}
\begin{align}\label{hopf}
 \Delta(\kad_{kl}) &= \sum_m \kad_{km}\mt \kad_{ml}\,,&  \counit(\kad_{kl}) &=   \delta_{kl}\,, & \antyp(\kad_{kl}) &= Ad^{\gotc}_{lk}\,;\\
\Delta(Ad^{\gotc}_{kl}) &=\sum_m Ad^{\gotc}_{km}\mt Ad^{\gotc}_{ml}\,,&  \counit(Ad^{\gotc}_{kl}) &= \delta_{kl}\,, &\antyp(Ad^{\gotc}_{kl}) &= \kad_{lk}\,.\nonumber
\end{align}
Let us define  $\displaystyle \tilde{A}_k:=\antyp(\kropka{c}_k)=-\sum_l Ad^{\gotc}_{lk} \kropka{c}_l\,$, then:  \notka{kappa-e}
\begin{align}\label{kappa-e}
\Delta(\tilde{A}_k)& =I\mt \tilde{A}_k + \sum_m \tilde{A}_m\mt Ad^{\gotc}_{mk}\,,&\,\antyp(\tilde{A}_k)&=-\sum_l \tilde{A}_l\, \kad_{kl}\,,& \,\kropka{c}_k&=- \sum_l \kad_{kl}\, \tilde{A}_l
\end{align}

Let us assume additionally, that $\gotg$ is equipped with invariant, non degenerate, symmetric bilinear form $k$ and that $k|_{\gotb}$ is non degenerate.
Thus we have two decompositions
\begin{equation*}
\gotg=\gotb\oplus \gotb^{\perp}=\gotb\oplus \gotc
\end{equation*}
%%%%%%%%%%%%%%%%%%%%
Let $P_{\gotb}$ be the projection on $\gotb$ defined by the first  decomposition and $P_{\gotc}$ projection on $\gotc$ defined by the second one. The next lemma is straightforward.
\notka{lema-reps}
\begin{lem}\label{lema-reps} (a) The restriction of $k$ to $\gotb^{\perp}$ is an isomorphism of $\gotb^{\perp}$  and $\gotb^0$.
For any $\kropka{c}\in\gotc$ and any $\kropka{f} \in \gotb^{\perp}$:  
$\displaystyle P_{\gotc}(I-P_{\gotb}) \kropka{c}=\kropka{c}\,,\,\,(I-P_{\gotb}) P_{\gotc}\kropka{f}=\kropka{f},$
in other words the restriction of $P_{\gotc}$ to $\gotb^{\perp}$ is an isomorphism of $\gotb^{\perp}$ and $\gotc$,  and the inverse mapping is the restriction of $I-P_{\gotb}$ to $\gotc$.\\
(b) The mapping  $\phi:=k \cdot (I-P_{\gotb})|_{\gotc} : \gotc\rightarrow \gotb^0$ is an isomorphism with the inverse $\phi^{-1}=P_{\gotc}\cdot( k^{-1}|_{\gotb^0}) $, moreover
for any $b\in B$ 
$$\phi\cdot  Ad^{\gotc}(b)\cdot \phi^{-1} =\kad(b)|_{\gotb^0}$$
(c) Let $(e_k)$ be o.n. basis in $\gotb^{\perp}$; it defines bases $(k(e_k))$ and $(P_{\gotc}(e_k))$ in $\gotb^0$ and $\gotc$, respectively. Then $\phi(P_{\gotc}(e_k))=k(e_k)$ and, consequently,
matrix elements of $Ad(b)|_{\gotb^{\perp}}$, $\kad(b)$ and $Ad^{\gotc}(b)$  in bases $(e_k)$, $(k(e_k))$ and $(P_{\gotc}(e_k))$, respectively, are equal.
\end{lem}
\dowl
%%%%%%%%%%%%%%%%%%%%%%%%%%%%%%%%%%%%%%%%%%%%%%%%%%%%%%%%%%%%%%%%%%%%%%%%%%%%%%%%%%%%%%%%%%%%%%%%%%%%%%%%%%%%%%%%%%

\end{document}